\documentclass[12pt]{article}

\usepackage{amsfonts,amsmath,latexsym,amssymb,txfonts,bm} 
\usepackage[american]{babel}
\makeatletter
\renewcommand{\@evenhead}{\raisebox{0pt}[\headheight][0pt]{\vbox{\hbox
to \textwidth{\thepage\hfil\strut\textit{\leftmark}}\hrule}}}
\renewcommand{\@oddhead}{\raisebox{0pt}[\headheight][0pt]{\vbox{\hbox
to \textwidth{\textit{\rightmark}\hfil\strut\thepage}\hrule}}}
\makeatother

\def\II{{\mathbb I}} 
\def\RR{{\mathbb R}} 
\def\CC{{\mathbb C}}

\def\tr{\mathrm{ tr\,}} 
\def\Tr{\mathrm{ Tr\,}}

\def\Ind{\mathrm{Ind\,}}
 
\def\vol{\mathrm{ vol\,}}

\def\diag{\mathrm{diag\,}} 
\def\Spin{\mathrm{Spin}}

\def\be{\begin{equation}} 
\def\ee{\end{equation}} 
\def\bea{\begin{eqnarray}} 
\def\eea{\end{eqnarray}} 
\def\bed{\begin{definition}{\ }}
\def\eed{\end{definition}}
\def\bd{\begin{description}}
\def\ed{\end{description}}
\def\bc{\begin{center}}
\def\ec{\end{center}}

\newtheorem{theorem}{Theorem}

\newtheorem{proposition}{Proposition}

\newtheorem{definition}{Definition}

\def\sideremark#1{\ifvmode\leavevmode\fi\vadjust{\vbox to0pt{\vss
\hbox to 0pt{\hskip\hsize\hskip1em
\vbox{\hsize2cm\tiny\raggedright\pretolerance10000
\noindent #1\hfill}\hss}\vbox to8pt{\vfil}\vss}}}

\begin{document}

\begin{titlepage}
\thispagestyle{empty}
\null
\hspace*{50truemm}{\hrulefill}\par\vskip-4truemm\par
\hspace*{50truemm}{\hrulefill}\par\vskip5mm\par
\hspace*{50truemm}{{\large\sc New Mexico Tech {\rm 
(\today)}}}\vskip4mm\par
\hspace*{50truemm}{\hrulefill}\par\vskip-4truemm\par
\hspace*{50truemm}{\hrulefill}
\par
\bigskip
\bigskip
%\par
%\hspace*{50truemm}{\LARGE\textbf{\textsf{DRAFT}}}
%\par
\vspace{3cm}
\centerline{\huge\bf Heat Kernel}
\bigskip
\centerline{\huge\bf on Homogeneous Bundles}
\bigskip
%\bigskip
\centerline{\huge\bf over Symmetric Spaces}
%\bigskip
\bigskip
\centerline{\Large\bf Ivan G. Avramidi}
\bigskip
\centerline{\it Department of Mathematics}
\centerline{\it New Mexico Institute of Mining and Technology}
\centerline{\it Socorro, NM 87801, USA}
\centerline{\it E-mail: iavramid@nmt.edu}
\bigskip
%\centerline{Revised on }
\medskip
%\maketitle 
\vfill

{\narrower
\par
% Abstract 

We consider Laplacians acting on sections of homogeneous vector bundles
over symmetric spaces.  By using an integral representation of the heat
semi-group  we find a formal solution for the heat kernel diagonal that
gives a generating function for the whole sequence of heat invariants.
We show explicitly that the obtained result correctly reproduces  the
first non-trivial heat kernel coefficient as well as the exact heat
kernel diagonals on two-dimensional sphere $S^2$ and the hyperbolic
plane $H^2$. We argue that the obtained formal solution correctly
reproduces the exact heat kernel diagonal after a suitable
regularization and analytical continuation.

\par}
\vfill

%{\vbox{
%\hrule
%\vspace{3pt}
%\hfil
%{\scriptsize\textit{
%\hfill\hfill\jobname.tex; \today; \timenow; p. \theoutputpage}}
%\hfil
%}}

\end{titlepage}

%=================================================================

\section{Introduction}
\setcounter{equation}0

The heat kernel is one of the most powerful tools in mathematical
physics and  geometric analysis (see, for example the books 
\cite{gilkey95,berline92,hurt83,avramidi00,kirsten01} and reviews 
\cite{avramidi87,camporesi90,avramidi99,avramidi02,vassilevich03}).
The  short-time  asymptotic expansion of the trace of the heat kernel
determines the spectral asymptotics of the differential operator. The
coefficients of this asymptotic expansion, called the  heat invariants,
are extensively used in geometric analysis, in particular, in spectral
geometry and index theorems proofs
\cite{gilkey95,berline92}. 

There has been a tremendous progress in the explicit calculation  of
spectral asymptotics in the last thirty years 
\cite{gilkey75,avramidi87,avramidi89,avramidi90,avramidi91,vandeven98,
yajima04}.  It seems that further progress in the study of spectral
asymptotics can be only achieved by restricting oneself to operators and
manifolds with high level of symmetry, in particular, homogeneous
spaces, which enables one to employ powerful algebraic methods.  In some
very special  particular cases, such as group manifolds, spheres,
rank-one symmetric spaces and split-rank symmetric spaces,  it is
possible to determine the spectrum of the Laplacian exactly and to
obtain closed formulas for the  heat kernel in terms of the root vectors
and their multiplicities 
\cite{anderson90,camporesi90,dowker70,dowker71,hurt83,fegan83}. The
complexity of the method crucially depends on the global structure of
the symmetric space, most importantly its rank. Most of the results for
symmetric spaces are obtained for rank-one symmetric spaces only
\cite{camporesi90}.

It is well known that heat invariants are  determined essentially by
local geometry. They are polynomial invariants in the curvature with
universal constants that do not depend on the global properties of the
manifold \cite{gilkey95}. It is this universal structure that we are
interested in this paper.   Our goal is to compute the heat kernel
asymptotics of the Laplacian acting on homogeneous vector bundles over
symmetric spaces. Related problems in a more general context are
discussed in \cite{avramidi94a,avramidi95,avramidi98}.

%=================================================================
%==================================================================
\section{Geometry of Symmetric Spaces}
\setcounter{equation}{0}

\subsection{Twisted Spin-Tensor Bundles}

In this section we introduce the basic concepts and fix notation.
Let  $(M,g)$ be an $n$-dimensional Riemannian manifold without boundary.
We assume that it is complete simply connected orientable and spin. We
denote the local coordinates on $M$ by $x^\mu$, with Greek indices
running over $1,\dots, n$. Let $e_{a}{}^\mu$ be a local orthonormal
frame defining a basis for the tangent space $T_xM$ so that
\be
g^{\mu\nu}=\delta^{ab}e_{a}{}^\mu e_{b}{}^\nu\,,
\ee
We denote the frame indices by low case Latin indices from the beginning
of the alphabet, which also run over $1,\dots,n$. The frame indices are
raised and lowered by the metric $\delta_{ab}$.  Let $e^{a}{}_\mu$ be the
matrix inverse to $e_{a}{}^\mu$, defining the dual basis in the
cotangent space $T_x^*M$, so that,  
\be
g_{\mu\nu}=\delta_{ab}e^{a}{}_\mu e^{b}{}_\nu\,.
\ee
The Riemannian volume element is defined as usual by
\be
d\vol=dx\,|g|^{1/2}\,,
\ee
where
\be
|g|=\det g_{\mu\nu}=(\det e_a{}^\mu)^2\,.
\ee
The spin connection $\omega^{ab}{}_{\mu}$ is defined in terms of
the orthonormal frame by 
\bea
\omega^{ab}{}_{\mu}
&=&
e^{a\mu}e^b{}_{\mu;\nu}=-e^{a}{}_{\mu;\nu} e^{b\mu}
\nonumber\\
&=&
e^{a\nu}\partial_{[\mu} e^{b}{}_{\nu]}
-e^{b\nu}\partial_{[\mu} e^{a}{}_{\nu]}
+e_{c\mu}e^{a\nu}e^{b\lambda}\partial_{[\lambda} e^{c}{}_{\nu]}\,,
\label{25xxx}
\eea
where the semicolon denotes the usual Riemannian covariant derivative
with the Levi-Civita connection.
The curvature of the spin connection is
\be
R^{a}{}_{b\mu\nu}=\partial_\mu\omega^{a}{}_{b\nu}
-\partial_\nu\omega^{a}{}_{b\mu}
+\omega^{a}{}_{c\mu}\omega^{c}{}_{b\nu}
-\omega^{a}{}_{c\nu}\omega^{c}{}_{b\mu}\,.
\ee
The Ricci tensor and the scalar curvature are defined by
\be
R_{\alpha\nu}= e_{a}{}^\mu e^{b}{}_\alpha R^{a}{}_{b\mu\nu}\,,
\qquad
R=g^{\mu\nu}R_{\mu\nu}=e_{a}{}^\mu e_{b}{}^\nu R^{ab}{}_{\mu\nu}\,.
\ee

Let $\mathcal{T}$ be a spin-tensor bundle realizing a representation
$\Sigma$  of the spin group \mbox{$\Spin(n)$}, the double covering of
the group $SO(n)$, with the fiber $\Lambda$. 
Let $\Sigma_{ab}$ be the generators of the 
orthogonal algebra ${\cal SO}(n)$, the Lie algebra of the orthogonal
group $SO(n)$, satisfying the following commutation relations
\be
[\Sigma_{ab},\Sigma_{cd}]
= -\delta_{ac}\Sigma_{bd}
+\delta_{bc}\Sigma_{ad}
+\delta_{ad}\Sigma_{bc}
-\delta_{bd}\Sigma_{ac}\,.
\label{280}
\ee

The spin connection induces a connection  on the bundle ${\cal T}$
defining the covariant derivative of smooth sections  $\varphi$ of the
bundle ${\cal T}$ by
\be
\nabla_\mu\varphi
=\left(\partial_\mu
+\frac{1}{2}\omega^{ab}{}_{\mu}\Sigma_{ab}\right)\varphi\,.
\ee
The commutator of covariant derivatives
defines the curvature of this connection via
\be
[\nabla_\mu,\nabla_\nu]\varphi
=\frac{1}{2}R^{ab}{}_{\mu\nu}\Sigma_{ab}\varphi\,.
\ee

As usual, the orthonormal frame,
$e^a{}_\mu$ and $e_a{}^\mu$, 
will be used to transform the coordinate
(Greek) indices to the orthonormal (Latin) indices. 
The covariant derivative  along the frame vectors is 
defined by
$\nabla_a=e_a{}^\mu\nabla_\mu$.
For example, with our notation, $\nabla_a\nabla_b T_{cd}
=e_a{}^\mu e_b{}^\nu e_c{}^\alpha e_d{}^\beta\nabla_\mu\nabla_\nu 
T_{\alpha\beta}$.

The metric $\delta_{ab}$ induces a positive definite fiber metric on
tensor bundles.  For Dirac spinors, the fiber metric is defined as
follows. First, one defines the Dirac matrices, $\gamma_{a}$, as
generators   of the Clifford algebra, (represented by $2^{[n/2]}\times
2^{[n/2]}$ complex matrices),
\be
\gamma_{a}\gamma_{b}+\gamma_{b}\gamma_{a}=2\delta_{ab}\II_S\,,
\ee
where $\II_S$ is the identity matrix in the spinor representation.
Then one defines the anti-symmetrized products of Dirac 
matrices
\be
\gamma_{a_1\dots a_k}=\gamma_{[a_1}\cdots\gamma_{a_k]}\,.
\ee
Then the matrices 
\be
\Sigma_{ab}=\frac{1}{2}\gamma_{ab}\,
\label{213mm}
\ee
are the generators of the orthogonal algebra ${\cal SO}(n)$
in the spinor representation.
The Hermitian conjugation of Dirac matrices defines a 
Hermitian matrix $\beta$
\footnote{The Dirac matrices $\gamma_{ab}$ and the 
spinor metric $\beta$ should not be confused with the
matrices $\gamma_{AB}$ and $\beta_{ij}$ defined below.}
by
\be
\gamma^{\dagger}_{a}=\beta\gamma_{a}\beta^{-1}\,,
\ee
which defines a Hermitian inner product
$\bar\psi\varphi=\psi^\dagger\beta\varphi$ in the vector space of
spinors. 
We also find the following important relation
\be
R^{ab}{}_{cd}\gamma_{ab}\gamma^{cd}
=-2 R^{ab}{}_{ab}\II_S
=-2R\;\II_S\,,
\label{117}
\ee
where $R$ is the scalar curvature.

In the present paper we will further  assume that $M$ is a {\it locally
symmetric space} with a Riemannian metric with the parallel curvature
\be
\nabla_\mu R_{\alpha\beta\gamma\delta}=0\,,
\ee
which means, in particular, that the curvature satisfies the 
integrability constraints
\be
R^{fg}{}_{ e a}R^e{}_{b cd} 
-R^{fg}{}_{ e b}R^e{}_{a cd} 
+ R^{fg}{}_{ ec}R^e{}_{d ab}
-R^{fg}{}_{ ed}R^e{}_{c ab}
= 0\,.
\label{312}
\ee

%========================================

Let $G_{YM}$ be a  compact Lie group (called a gauge group). It
naturally defines the principal fiber bundle over the manifold $M$  with
the structure group $G_{YM}$. We consider a representation of the
structure group $G_{YM}$ and the associated vector bundle through
this representation with the same structure group $G_{YM}$ whose typical
fiber is a $k$-dimensional  vector space $W$. Then for any spin-tensor
bundle ${\cal T}$ we define the twisted spin-tensor bundle $\mathcal{V}$
via the twisted product of the bundles ${\cal W}$ and ${\cal T}$. The
fiber of the bundle ${\cal V}$ is $V=\Lambda\otimes W$ so that the
sections of the bundle ${\cal V}$ are represented locally by $k$-tuples
of spin-tensors.

Let ${\cal A}$ be a connection one form on the bundle ${\cal W}$
(called Yang-Mills or gauge connection) taking values in the Lie
algebra ${\cal G}_{YM}$ of the gauge group $G_{YM}$. Then the total
connection on the bundle ${\cal V}$ is defined by
\be
\nabla_\mu\varphi
=\left(\partial_\mu
+\frac{1}{2}\omega^{ab}{}_\mu \Sigma_{ab}\otimes \II_W
+\II_\Lambda\otimes{\cal A}_\mu\right)\varphi\,,
\ee
and the total curvature $\Omega$ of the bundle ${\cal V}$ is defined by
\be
[\nabla_\mu,\nabla_\nu]\varphi
=\Omega_{\mu\nu}\varphi\,,
\label{219mm}
\ee
where
\be
\Omega_{\mu\nu}=
\frac{1}{2}R^{ab}{}_{\mu\nu} \Sigma_{ab}
+{\cal F}_{\mu\nu}\,,
\label{220}
\ee
and
\be
{\cal F}_{\mu\nu}
=\partial_\mu {\cal A}{}_\nu
-\partial_\nu {\cal A}{}_\mu
+[{\cal A}_\mu,{\cal A}_\mu]
\ee
is the curvature of the Yang-Mills connection.

We also consider the bundle of endomorphisms of the bundle ${\cal V}$.
The covariant derivative  of sections
of this bundle is defined by
\be
\nabla_\mu X
=\left(\partial_\mu
+\frac{1}{2}\omega^{ab}{}_\mu \Sigma_{ab}
\right)X+[{\cal A}_\mu,X]\,,
\ee
and the commutator of covariant derivatives is equal to
\be
[\nabla_\mu,\nabla_\nu]X
=\frac{1}{2}R^{ab}{}_{\mu\nu} \Sigma_{ab}X+[{\cal F}_{\mu\nu},X]\,.
\label{225}
\ee

In the following we will consider {\it homogeneous vector bundles} with
parallel bundle curvature
\be
\nabla_\mu {\cal F}_{\alpha\beta}=0\,,
\ee
which means that the curvature satisfies the integrability constraints
\be 
[{\cal F}_{cd},{\cal F}_{ab}]  - R^f{}_{acd}{\cal F}_{fb} -
R^f{}_{bcd}{\cal F}_{af} = 0\,. 
\label{227} 
\ee

%=================================================================
%\section{Symmetric Spaces}
%\setcounter{equation}0

%===================================================================
\subsection{Normal Coordinates}

Let $x'$ be a fixed point in $M$ and ${\cal U}$ be a sufficiently small
coordinate patch containing the point $x'$. Then every point $x$ in
${\cal U}$ can be connected with the point $x'$ by a unique geodesic. We
extend the local orthonormal frame $e_a{}^\mu(x')$ at the point $x'$ to
a local  orthonormal frame $e_a{}^\mu(x)$ at the point $x$ by parallel
transport
\be
e_a{}^\mu(x)=g^\mu{}_{\nu'}(x,x')e_a{}^{\nu'}(x')\,,
\ee
\be
e^a{}_\mu(x)=g_\mu{}^{\nu'}(x,x')e^a{}_{\nu'}(x')\,,
\ee 
where $g^\mu{}_{\nu'}(x,x')$  is the operator of parallel transport of
vectors along the geodesic from the point $x'$ to the point $x$. Of
course, the frame $e_a{}^\mu$ depends on the fixed point $x'$ as a
parameter. Here and everywhere below the coordinate indices of the
tangent space at the point $x'$ are denoted by primed Greek letters. 
They are raised and lowered by the metric tensor $g_{\mu'\nu'}(x')$ at
the point $x'$.  The derivatives with respect to $x'$ will be denoted by
primed Greek indices as well. 

The parameters of the geodesic connecting the points $x$ and $x'$,
namely the unit tangent vector at the point $x'$  and the length of the
geodesic, (or, equivalently, the tangent vector at the point $x'$ with
the norm equal to the length of the geodesic), provide normal coordinate
system for ${\cal U}$. Let $d(x,x')$ be the geodesic distance between
the points $x$ and $x'$ and $\sigma(x,x')$ be a two-point function 
defined by
\be
\sigma(x,x')=\frac{1}{2}[d(x,x')]^2\,.
\ee
Then the derivatives $\sigma_{;\mu}(x,x')$ 
and $\sigma_{;\nu'}(x,x')$ are the tangent vectors to
the geodesic connecting the points $x$ and $x'$ at the points $x$
and $x'$ respectively pointing in opposite directions; one is obtained 
from another by parallel transport
\be
\sigma_{;\mu}=-g_\mu{}^{\nu'}\sigma_{;\nu'}\,.
\ee
Here and everywhere below the semicolon denotes the covariant
derivative.

The operator of parallel transport satisfies the equation
\be
\sigma^{;\mu}\nabla_{\mu} g^\alpha{}_{\beta'}=0\,,
\ee
with the initial conditions
\be
g^\alpha{}_{\beta'}\Big|_{x=x'}=\delta^\alpha_\beta\,.
\ee
It can be expressed in terms of the
local parallel frame 
\be
g^\mu{}_{\nu'}(x,x')=e_a{}^\mu(x) e^a{}_{\nu'}(x')\,,
\ee
\be
g_\mu{}^{\nu'}(x,x')=e^a{}_\mu(x)e_a{}^{\nu'}(x')\,.
\ee

Now, let us define the quantities
\be
y^a
=e^{a}{}_{\mu}\sigma^{;\mu}
=-e^{a}{}_{\mu'}\sigma^{;\mu'}\,,
\ee
so that
\be
\sigma^{;\mu}=e_{a}{}^{\mu}y^a
\qquad
\mbox{and}
\qquad
\sigma^{;\mu'}=-e_{a}{}^{\mu'}y^a
\,.
\ee
Notice that $y^a=0$ at $x=x'$.
Further, we have
\be
\frac{\partial y^a}{\partial x^\nu}
=-e^{a\,\mu'}\sigma_{;\nu\mu'}\,,
\ee
so that
the Jacobian of the change of variables
is
\be
\det\left(\frac{\partial y^a}{\partial x^\nu}\right)
= |g|^{-1/2}(x')\det[-\sigma_{;\nu\mu'}(x,x')]\,.
\ee

The geometric parameters $y^a$ are nothing but the normal coordinates.
By using the Van Vleck-Morette determinant 
defined by 
\footnote{Do not confuse it with the Laplacian $\Delta$ defined below.}
\be
\Delta(x,x')=|g|^{-1/2}(x')|g|^{-1/2}(x)
\det[-\sigma_{;\nu\mu'}(x,x')]\,,
\ee
we can write the
Riemannian volume element in the form
\be
d\vol
=dy\; \Delta^{-1}(x,x')\,.
\ee

Let ${\cal P}(x,x')$ be the operator of parallel transport
of sections of the bundle ${\cal V}$ from the point $x'$ to the point $x$.
It satisfies the equation
\be
\sigma^{;\mu}\nabla_\mu {\cal P}=0\,,
\ee
with the initial condition
\be
{\cal P}\Big|_{x=x'}=\II_{V}\,.
\ee

Any spin-tensor $\varphi$ can be now expanded
in the covariant Taylor series
\be
\varphi(x)=
{\cal P}(x,x')\sum_{k=0}^\infty \frac{1}{k!}
\left[
\nabla_{(c_1}\cdots\nabla_{c_k)}
\varphi\right](x')
y^{c_1}\cdots y^{c_k}
\,.
\label{243}
\ee
Therefrom it is clear, in particular, that the frame components of a 
parallel spin-tensor are simply constant.

In symmetric spaces
one can compute the Van Vleck-Morette determinant
explicitly in terms of the curvature.
Let  $K$ be a $n\times n$ matrix with the entries  
\be
K^a{}_b=R^a{}_{cbd}y^cy^d\,.
\ee
Then \cite{avramidi91,avramidi87,avramidi00}
\be
\frac{\partial y^a}{\partial x^\nu}
=\left(\frac{\sqrt{K}}{\sin\sqrt{K}}\right)^a{}_b\; e^b{}_\nu\,,
\label{244mm}
\ee
and, therefore,
\be
\Delta(x,x')=
\det{}_{TM}\left(\frac{\sqrt{K}}{\sin\sqrt{K}}\right)\,.
\label{245xx}
\ee
Thus, the Riemannian volume element in symmetric spaces takes the following 
form
\be
d\vol=dy\;\det{}_{TM}\left(\frac{\sin\sqrt{K}}{\sqrt{K}}\right)\,.
\ee

%=============================================

The matrix $(\sin\sqrt{K})/\sqrt{K}$ determines the orthonormal frame in
normal coordinates, and the square of this matrix determines the metric
tensor in normal coordinates,
\be
ds^2=\left(\frac{\sin^2\sqrt{K}}{K}\right)_{ab} dy^a\, dy^b\,.
\ee

Let us define an endo-morphism valued $1$-form $\tilde{\cal A}_a$ by
the equation
\be
\nabla_\nu {\cal P}=
{\cal P}\tilde {\cal A}_a e^a{}_{\mu'}\sigma^{;\mu'}{}_{\nu}
\,.
\ee
Then for bundles with parallel curvature over symmetric spaces one can find
it explicitly
\cite{avramidi91,avramidi87,avramidi00}
\be
\tilde{\cal A}_{a}=-{\cal F}_{bc}y^c 
\left(\frac{\II-\cos\sqrt{K}}{K}\right)^b{}_a\,.
\ee
This object determines the gauge connection in normal coordinates,
\be
{\cal A}=-{\cal F}_{bc}y^c 
\left(\frac{\II-\cos\sqrt{K}}{K}\right)^b{}_a\, dy^a\,.
\ee

This means that all connections on a homogeneous bundle are essentially
the same. In particular, the spin connection one-form in normal
coordinates has the form
\be
\omega^a{}_b=-R^a{}_{bcd}y^d
\left(\frac{\II-\cos\sqrt{K}}{K}\right)^c{}_e dy^e\,.
\label{251xxx}
\ee

%==================================================================

{\bf Remarks.}
Two remarks are in order here. First, strictly speaking, normal
coordinates can be only defined locally, in geodesic balls of radius
less than the injectivity radius of the manifold. However, for symmetric
spaces normal coordinates cover the whole manifold except for a set of
measure zero where they become singular \cite{camporesi90}. This set is
precisely the set of points conjugate to the fixed point $x'$ (where 
$\Delta^{-1}(x,x')=0$) and of points that can be connected to the point
$x'$ by multiple geodesics. In any case, this set is a set of measure
zero and, as we will show below, it can be dealt with  by some
regularization technique. Thus, we will use the normal coordinates
defined above for the whole manifold. Second, for compact manifolds (or
for manifolds with compact submanifolds) the range of some normal
coordinates is also compact, so that if one allows them to range over
the whole real line $\RR$, then the corresponding compact submanifolds
will be covered infinitely many times.

%==================================================================
\subsection{Curvature Group of a Symmetric Space}

We assumed that the manifold $M$ is locally symmetric. Since we also
assume that it is simply
connected and complete, it is a globally symmetric space (or simply
symmetric space) \cite{wolf72}.  A symmetric space is said to be
compact, non-compact or Euclidean if all sectional curvatures are
positive, negative or zero. A generic symmetric space has the structure
\be
M=M_0\times M_s\,,
\ee
where $M_0=\RR^{n_0}$ and
 $M_s$ is a semi-simple symmetric space; it is a product
of a compact symmetric space $M_+$ and a non-sompact symmetric space $M_-$,
\be
M_s=M_+\times M_-\,.
\ee
Of course, the dimensions must satisfy the relation  $n_0+n_s=n$,
where $n_s=\dim M_s$.

Let $\Lambda_2$ be the vector space of $2$-forms on $M$
at a fixed point $x'$. It has the
dimension $\dim \Lambda_2=n(n-1)/2$, and the inner product in
$\Lambda_2$ is defined by 
\be
\left<X,Y\right>=\frac{1}{2}X_{ab}Y^{ab}\,.
\ee
The Riemann curvature tensor naturally defines the curvature operator 
\be
{\rm Riem}: \Lambda_2\to \Lambda_2
\ee
by
\be
({\rm Riem}\, X)_{ab}=\frac{1}{2}R_{ab}{}^{cd}X_{cd}\,.
\ee
This operator is symmetric and has real eigenvalues which determine the
principal sectional curvatures. Now, let ${\rm Ker}\,({\rm Riem})$
and ${\rm Im}\,({\rm Riem})$ be the
kernel and the range 
of this operator and 
\be
p=\dim {\rm Im}({\rm Riem})
=\frac{n(n-1)}{2}-\dim {\rm Ker}\,({\rm Riem})\,\,.
\ee
Further, let $\lambda_i$, $(i=1,\dots, p)$, be the non-zero eigenvalues,
and $E^i{}_{ab}$ be the corresponding orthonormal eigen-two-forms.
Then the components of the curvature tensor
can be presented in the form \cite{avramidi96}
\be
R_{abcd} = \beta_{ik}E^i{}_{ab}E^k{}_{cd}\,,
\label{236}
\ee
where  $\beta_{ik}$ is a symmetric, in fact, diagonal, nondegenerate  
$p\times p$ matrix
\be
(\beta_{ik})=\diag(\lambda_1,\dots,\lambda_p)\,.
\ee
Of course, the zero eigenvalues of the curvature operator correspond to
the flat subspace $M_0$, the positive ones correspond to the compact
submanifold $M_+$ and the negative ones to the non-compact submanifold
$M_-$. Therefore, ${\rm Im}\,({\rm Riem})=T_x M_s$.

In the following the
Latin indices from the middle of the alphabet will be used to denote
tensors in ${\rm Im}({\rm Riem})$; 
they should not be confused with the Latin indices from
the beginning of the alphabet which denote tensors in $M$.
They will be raised and lowered with
the matrix $\beta_{ik}$ and its inverse
\be
(\beta^{ik})=\diag(\lambda_1^{-1},\dots,\lambda_p^{-1})
\,.
\ee

Next, we define the traceless $n\times n$
matrices $D_i=(D^a{}_{ib})$,
where
\be
D^a{}_{ib}=-\beta_{ik}E^k{}_{cb}\delta^{ca}\,.
\ee
Then
\be
R^a{}_{bcd}=-D^a{}_{ib}E^i{}_{cd}\,,
\ee
\be
R^a{}_b{}^c{}_d=\beta^{ik}
D^a{}_{ib}D^c{}_{kd}\,,
\ee
\be
R^a{}_b=-\beta^{ik}D^a{}_{ic}D^c{}_{kb},
\ee
\be
R=-\beta^{ik}D^a{}_{ic}D^c{}_{ka}\,.
\ee
Also, we have identically,
\be
D^a{}_{j[b} E^j{}_{cd]} = 0\,.
\label{248}
\ee

The matrices $D_i$ are known to be the generators of the  holonomy
algebra, ${\cal H}$, i.e. the Lie algebra of the restricted holonomy 
group, $H$,
\be
[D_i, D_k] = F^j{}_{ik} D_j\,,
\label{310}
\ee
where $F^j{}_{ik}$ are the structure constants of the holonomy group.
The structure constants of the holonomy group define the
$p\times p$ matrices $F_i$, by $(F_i)^j{}_k=F^j{}_{ik}$, which generate
the adjoint representation  of the holonomy algebra,
\be
[F_i, F_k] = F^j{}_{ik} F_j\,.
\ee
These commutation relations follow directly from the Jacobi identities
\be
F^I{}_{j[k}F^j{}_{ml]}=0\,.
\label{252}
\ee

For symmetric spaces the introduced quantities satisfy additional
algebraic constraints. The most important consequence of the eq.
(\ref{312}) is the equation
%\footnote{We correct here a sign misprint in eq. (3.13) in 
%\cite{avramidi96}.}
\cite{avramidi96}
\be
E^i{}_{a c} D^c{}_{kb}
-E^i{}_{b c} D^c{}_{ka} 
= F^i{}_{kj}E^j{}_{ab}\,.
\label{313}
\ee
It is this equation that makes a generic Riemannian manifold
a symmetric
space.

Now, by using the eqs. (\ref{310}) and (\ref{313}) one can prove the 
following:
\begin{proposition}
The matrix $\beta_{ik}$ is $H$-invariant and satisfies the equation
\be
\beta_{ik}F^k{}_{jl}+\beta_{lk}F^k{}_{ji}=0\,.
\label{270}
\ee
\end{proposition}
This means that the matrices $F_i$ satisfy the transposition rule
\be
(F_i)^T=-\beta F_i \beta^{-1}\,,
\ee
which simply means that 
the adjoint and the coadjoint representations of
the holonomy algebra ${\cal H}$ are equivalent.
In particular, this means that the matrices $F_i$ are
traceless.
Such an algebra is called {\it compact}
\cite{barut77}.

Another consequence of the eq. (\ref{313}) are the
identities
\be
D^{a}{}_{i[b} R_{c]ade}+D^{a}{}_{i[d} R_{e]abc}= 0,
\label{257}
\ee
\be
R^a{}_{c}D^c{}_{ib}=D^a{}_{ic}R^c{}_{b}\,.
\ee
This means, in particular, that the Ricci tensor matrix commutes with all 
matrices $D_i$ and is, therefore, an invariant matrix of the holonomy 
algebra. Thus,
\be
R^a{}_b =\frac{1}{n_s}h^a{}_b R\,,
\ee
where $h^a{}_b$ is a projection (a symmetric idempotent parallel
tensor) to
the subspace $T_x M_s$ of the tangent space of dimension $n_s$, that is,
\be
h_{ab}=h_{ba}\,,\qquad
h^a{}_b h^b{}_c=h^a{}_c\,,\qquad
h^a{}_a=n_s\,.
\ee
It is easy to see that the tensor $h_{ab}$ is nothing but the metric tensor
on the semi-simple subspace $T_xM_s$.

Since the curvature exists only in the semi-simple submanifold $M_s$, the
components of the curvature tensor $R_{abcd}$, as well as the tensors
$E^i{}_{ab}$,  are non-zero only in the semi-simple subspace $T_xM_s$.
Let
\be
q^a{}_{b}=\delta^a{}_b-h^a{}_b\,
\ee
be the projection tensor to the flat subspace $\RR^{n_0}$ such that
\be
q_{ab}=q_{ba}\,, \qquad q^a{}_b q^b{}_c=q^a{}_c\,,\qquad
q^a{}_a=n_0\,, \qquad
q^a{}_b h^b{}_c=0\,.
\ee
Then 
\be
R_{abcd}q^a{}_e=R_{ab}q^a{}_e=E^i{}_{ab}q^a{}_e=D^a{}_{ib}q^b{}_e
=D^a{}_{ib} q_a{}^e=0\,.
\ee 

Now, we introduce a new type of indices, the capital Latin indices,
$A,B,C,\dots,$ which split according to $A=(a,i)$ and run from $1$ to
$N=p+n$. We define new quantities $C^A{}_{BC}$ by
\be
C^i{}_{ab}=E^i{}_{ab}, \qquad 
C^a{}_{ib}=-C^a{}_{bi}=D^a{}_{ib}, \qquad 
C^i{}_{kl}=F^i{}_{kl}\,,
\label{317}
\ee
all other components being zero.
Let us also introduce rectangular $p\times n$ matrices 
$T_a$ by $(T_a)^j{}_c=E^j{}_{ac}$ and the $n\times p$ matrices
$\bar T_a$ by $(\bar T_a)^b{}_i=-D^b{}_{ia}$. 
Then we can define $N\times N$ matrices $C_A=(C_a,C_i)$
\be
C_a = \left(
\begin{array}{cc}
0 & \bar T_a \\
T_a & 0 \\
\end{array}
\right)\,,
\qquad
C_i = \left(
\begin{array}{cc}
D_i & 0 \\             
0 & F_i\\
\end{array}
\right),
\label{318}
\ee
so that $(C_A)^B{}_C=C^B{}_{AC}$. 

\begin{theorem}
The quantities $C^A{}_{BC}$ satisfy the Jacobi identities
\be
C^A{}_{B[C}C^C{}_{DE]}=0\,.
\ee
This means that the matrices $C_A$ satisfy the commutation
relations
\be
[C_A, C_B]=C^C{}_{AB}C_C\,,
\label{320}     
\ee
or, in more details,
\bea
[C_a,C_b]&=&E^i{}_{a b}C_i,
\\{}
[C_i,C_a]&=&D^b{}_{ia}C_b,
\\{}
[C_i, C_k] &=& F^j{}_{i k} C_j\,,
\label{321}
\eea
and generate the adjoint representation of a Lie algebra ${\cal G}$ 
with the structure constants $C^A{}_{BC}$.
\end{theorem}

\noindent{\bf Proof.} This can be proved by
using the eqs. (\ref{248}),  (\ref{310}),  (\ref{252})
and (\ref{313}) \cite{avramidi96}.

For the lack of a better name we call the algebra ${\cal G}$ the {\it
curvature algebra}. As it will be clear from the next section it is a
subalgebra of the total isometry algebra of the symmetric space. It
should be clear  that the holonomy algebra ${\cal H}$ is the subalgebra
of the curvature algebra  ${\cal G}$.  The curvature algebra exists only
in symmetric spaces; it is the eq. (\ref{313}) that closes this algebra.

Next, we define a symmetric nondegenerate
$N\times N$ matrix
\be
(\gamma_{AB}) = 
\left(
\begin{array}{cc}
\delta_{ab} & 0 \\
0 & \beta_{ik} \\
\end{array}
\right)
=\diag\left(\underbrace{1,\dots,1}_{n},\lambda_1,\dots,\lambda_p\right)
\,.
\label{319}
\ee
This matrix and its inverse 
$(\gamma^{AB})=
\left(
\begin{array}{cc}
\delta^{ab} & 0 \\
0 & \beta^{ik} \\
\end{array}
\right)
=
\diag(\underbrace{1,\dots,1}_{n},\lambda_1^{-1},\dots,\lambda_p^{-1})$
will be
used to lower and to raise the capital Latin indices. 

Finally, by using the eqs. (\ref{313}) and (\ref{270}) one can show 
the following:
\begin{proposition}
The matrix
$\gamma_{AB}$ is $G$-invariant and satisfies the equation
\be
\gamma_{AB}C^B{}_{CD}+\gamma_{DB}C^B{}_{CA}=0\,.
\label{288xx}
\ee
\end{proposition}
In matrix notation this equation takes the form
\be
(C_A)^T=-\gamma C_A\gamma^{-1}\,,
\label{289xx}
\ee
which means that the adjoint and the coadjoint representations
of the curvature group are equivalent. In particular,
the matrices $C_A$ are traceless.

Thus the curvature algebra ${\cal G}$ is compact; it is a direct sum
of two ideals, 
\be
{\cal G}={\cal G}_0\oplus {\cal G}_s,
\label{290xx}
\ee
an Abelian center ${\cal G}_0$ of dimension $n_0$ 
and a semi-simple algebra ${\cal G}_s$ of
dimension $p+n_s$.

It is worth mentioning that although the holonomy algebra ${\cal H}$ is
compact the (indefinite, in general) metric, $\beta_{ij}$, introduced
above is {\it not} equal to the (positive definite) Cartan-Killing form,
$\rho_{ij}$, defined by
\be
\tr_{TM} D_i D_k=D^a{}_{ib}D^b{}_{ka}=-\rho_{ik}\,,
\label{289a}
\ee
so that
\be
\rho_{ik}=\diag(\lambda_1^2,\dots,\lambda_p^2)\,,
\ee
and
\be
\beta^{ik}\rho_{ik}=R\,.
\ee
Similarly, the generators $F_i$ satisfy
\be
\tr_H F_i F_k=F^j{}_{im}F^m{}_{kj}=-4\frac{R_H}{R}\rho_{ik}\,,
\ee
where
\be
R_H=-\frac{1}{4}\beta^{ik}F^j{}_{im}F^m{}_{kj}\,.
\label{295vv}
\ee
The Killing-Cartan form $\tr_G C_A C_B$ for the curvature
algebra ${\cal G}$ is defined by
\be
\tr_G C_a C_b=-\frac{2}{n_s}h_{ab}R\,,
\ee
\be
\tr_G C_i C_j=-\left(1+4 \frac{R_H}{R}\right)\rho_{ij}\,,
\ee
\be
\tr_G C_a C_i=0\,.
\ee
Notice that it is degenerate and is not equal to the metric
$\gamma_{AB}$.

%===================================================================
\subsection{Killing Vectors Fields}

We will use extensively the isometries of the symmetric space $M$. 
We follow the approach developed in 
\cite{avramidi96,avramidi87,avramidi91,avramidi00}.
The
generators of isometries are the Killing vector fields
$\xi$ defined by the equation
\be
\nabla_\mu\xi{}_{\nu}+\nabla_\nu\xi{}_{\mu}=0\,.
\label{22a}
\ee
The integrability conditions for this equation are
\be
R_{\alpha\beta\mu[\lambda}\nabla_{\nu]}\xi^\mu
+R_{\lambda\nu\mu[\beta}\nabla_{\alpha]}\xi^\mu=0\,.
\label{276}
\ee
By differentiating this equation, commuting derivatives and using
curvature identities we obtain
\be
\nabla_\mu\nabla_\nu\xi{}^{\lambda}
=-R^{\lambda}{}_{\nu\alpha\mu}\xi{}^{\alpha}\,,
\label{22}
\ee
which means, in particular,
\be
\Delta\xi^\lambda=-R^\lambda{}_\alpha\xi^\alpha\,.
\ee
By induction we obtain
\bea
\nabla_{\mu_{2k}}\cdots\nabla_{\mu_{1}}\xi{}^\lambda
&=& (-1)^k 
R^\lambda{}_{\mu_1\alpha_1\mu_2}
R^{\alpha_1}{}_{\mu_3\alpha_2\mu_4}\cdots 
R^{\alpha_{k-1}}{}_{\mu_{2k-1}\alpha_{k}\mu_{2k}}
\xi{}^{\alpha_{k}}\,,
\label{250}
\\
\nabla_{\mu_{2k+1}}\cdots\nabla_{\mu_{1}}\xi{}^\lambda 
&=& (-1)^k 
R^\lambda{}_{\mu_1\alpha_1\mu_2}
R^{\alpha_1}{}_{\mu_3\alpha_2\mu_4}\cdots 
R^{\alpha_{k-1}}{}_{\mu_{2k-1}\alpha_{k}\mu_{2k}}\nabla_{\mu_{2k+1}}
\xi{}^{\alpha_{k}}\,.
\label{251}
%\nonumber\\
\eea
These derivatives determine all coefficients of the covariant Taylor
series (\ref{243})
for the Killing vectors, and therefore, every Killing vector in a
symmetric space has the form
\be
\xi^a(x)=\left(\cos\sqrt K\right)^a{}_b \xi^{b}(x')
+\left(\frac{\sin\sqrt{K}}{\sqrt{K}}\right)^a{}_b y^c
\xi^{b}{}_{;c}(x')
\,,
\ee
or
\be
\xi(x)=\left\{
\left(\sqrt{K}\cot\sqrt K\right)^a{}_b \xi^{b}(x')
+\xi^{a}{}_{;c}(x')y^c
\right\}
\frac{\partial}{\partial y^a}
\,.
\ee
Thus, Killing vector
fields at any point $x$ are determined by their values $\xi^a(x')$ and
the values of their derivatives $\xi^a{}_{;c}(x')$ at the fixed point
$x'$. 

Similarly we can obtain the derivatives of the Killing vectors
\bea
\xi^a{}_{;b}(x)
&=&
\xi^a{}_{;b}(x')
-R^a{}_{bcd}y^d\left(\frac{1-\cos\sqrt{K}}{K}\right)^c{}_ey^f \xi^e{}_{;f}(x')
\nonumber\\
&&
-R^a{}_{bcd}y^d \left(\frac{\sin\sqrt{K}}{\sqrt{K}}\right)^c{}_e \xi^e(x')
\,.
\label{2105}
\eea

The set of all Killing vector fields forms a representation of the 
isometry algebra, the Lie algebra of the isometry group
of the manifold $M$. We define two
subspaces of the isometry algebra. One subspace is formed by Killing
vectors satisfying the initial
conditions 
\be
\nabla_\mu\xi^\nu\Big|_{x=x'}=0\,,
\label{289}
\ee
and another subspace is formed by the Killing vectors 
satisfying the
initial conditions
\be
\xi^\nu\Big|_{x=x'}=0\,.
\label{290}
\ee
We will call the Killing vectors from the first subspace {\it 
translations} and the
Killing vectors from the second group 
{\it rotations}. However, this should not
be understood literally.

One can easily show that the initial values $\xi^a(x')$ are
independent and, therefore, there are $n$ such parameters.  Thus, there
are $n$ linearly independent translations, which can be chosen in the form
\be
P_a=\left(\sqrt{K}\cot\sqrt{K}\right)^b{}_a\frac{\partial}{\partial y^b}\,,
\label{226}
\ee
so that
\be
e^b{}_\mu P_a{}^\mu\big\vert_{x=x'}=\delta^b{}_a\,,\qquad
P_a{}^\mu{}_{;\nu}\big\vert_{x=x'} = 0,
\label{237}
\ee

It is worth poiting out that the nature of the lower index of the
Killing vectors $P_a{}^\mu$ is different from the frame indices.
This means, in particular, that the covariant derivative of $P_a{}^\mu$ 
does not include the spin connection associated 
with the lower index.
In other words, $P_a{}^\mu$ are just $n$ vectors and not the 
components of a $(1,1)$ tensor.

On the other hand, the initial values of the derivatives
$\xi^a{}_{;c}(x')$ are not independent because of the constraints
(\ref{276}). These constraints are valid only in the semi-simple 
subspace $T_xM_s$.  However, in this subspace, due to the identity
(\ref{257}), it should be clear that there are
$p$ linearly independent rotations
\be
L_i=-D^b{}_{ia}y^a\frac{\partial}{\partial y^b}\,,
\label{229}
\ee
satisfying the initial conditions
\be
L_i{}^\mu\Big|_{x=x'}=0\,,\qquad
e^a{}_\mu e_b{}^\nu  L_i{}^\mu{}_{;\nu}\Big|_{x=x'}
=-D^a{}_{ib}\,.
\label{238}
\ee

More generally, by using (\ref{2105}) we also obtain
\be
e^a{}_\mu e_b{}^\nu P_e{}^\mu{}_{;\nu}
=-R^a{}_{bcd}y^d\left(\frac{\sin\sqrt{K}}{\sqrt{K}}\right)^c{}_e\,,
\ee
\be
e^a{}_\mu e_b{}^\nu L_i{}^\mu{}_{;\nu}
=-D^a{}_{ib}+R^a{}_{bcd}y^d
\left(\frac{1-\cos\sqrt{K}}{K}\right)^c{}_e y^f D^e{}_{if}\,.
\ee
This means, in particular, that the derivatives of all 
Killing vectors have the form
\be
\xi_A{}^a{}_{;b}=-D^a{}_{ib}\eta_A{}^i\,,
\label{2115}
\ee
where $\eta_A{}^i$ are defined by
\be
\eta_A{}^i=\alpha^{ij}\xi_A{}^a{}_{;b}D^b{}_{ja}\,,
\label{2116}
\ee
and the matrix $\alpha^{ij}=(\rho_{ij})^{-1}$ is the inverse matrix
of the Cartan-Killing form $\rho$ defined by (\ref{289a}).
Notice that
\be
\eta_a{}^i\Big|_{x=x'}=0\,,
\qquad
\eta_j{}^i\Big|_{x=x'}=\delta^i_j\,.
\ee

Then, from the eq. (\ref{22}) we also immediately obtain
\be
\eta_A{}^i{}_{;b}=-E^i{}_{ab}\xi_A{}^a\,.
\label{2117}
\ee

By adding the trivial Killing vectors for flat subspaces we find that the
dimension of the rotation subspace is equal to
\be
p+n_0n_s+\frac{n_0(n_0-1)}{2}\,.
\ee
Here $n_0n_s$ is the number of mixed rotations between $M_0$ and $M_s$ 
and $n_0(n_0-1)/2$
is the number of rotations of $M_0$.
Since $p\le n_s(n_s-1)/2$, then the above number of rotations is
less or equal to $n(n-1)/2$ as it should be (recall that $n=n_0+n_s$).

In the following we will need only  the Killing vectors $P_a$ and $L_i$
defined above. We introduce the following notation $(\xi_A)=(P_a,L_i)$. 

\begin{theorem}
The Killing vector fields $\xi_A$ satisfy the commutation relations
\be
[\xi_A, \xi_B]=C^C{}_{AB}\xi_C\,,
\label{320a}
\ee
or, in more detail,
\bea
[P_a,P_b]&=&E^i{}_{a b}L_i,
\\{}
[L_i,P_a]&=&D^b{}_{ia}P_b,
\\{}
[L_i, L_k] &=& F^j{}_{i k} L_j\,.
\label{321a}
\eea
\end{theorem}

\noindent
\noindent\noindent{\bf Proof.} This can be proved by using the explicit
form of the Killing vector fields obtained above \cite{avramidi96}.

Notice that they {\it do not} generate  the complete isometry algebra of
the symmetric space $M$ but rather they form a representation of the
curvature algebra ${\cal G}$ introduced in the previous section, which
is a subalgebra of the total isometry algebra.

It is clear that the Killing vector fields $L_i$ form a
representation of the holonomy algebra ${\cal H}$, which is the isotropy
algebra of the semi-simple submanifold $M_s$, and a  subalgebra of the
total isotropy algebra of the symmetric space $M$.

\begin{proposition}
There holds
\be
\xi_A^c{}_{;a}\xi_B{}^b{}_{;c}
-\xi_B^c{}_{;a}\xi_A{}^b{}_{;c}
=C^C{}_{AB}\xi_C{}^b{}_{;a}-R^b{}_{acd}\xi_A{}^c\xi_B{}^d\,,
\label{2124x}
\ee
and
\be
F^j{}_{ik}\eta_A{}^i\eta_B{}^k
=C^C{}_{AB}\eta_C{}^j
-E^j{}_{cd}\xi_A^c\xi_B{}^d\,.
\label{2124}
\ee
\end{proposition}

\noindent
{\bf Proof.} By differentiating the eq. (\ref{320a}) and using
(\ref{22}) we obtain (\ref{2124x}). Finally, by using
(\ref{2115}) and the holonomy algebra  (\ref{310}) we obtain
(\ref{2124}).

%==================================
%\subsubsection{Bilinear Identities}

Now, we derive some bilinear identities that we will need in the present
paper.

\begin{proposition}
The Killing vector fields satisfy the equation
\be
\gamma^{AB}\xi_A{}^\mu\xi_B{}^\nu
=\delta^{ab} P_a{}^\mu P_b{}^\nu
+\beta^{ik}L_i{}^\mu L_k{}^\nu
=g^{\mu\nu}\,,
\label{267}
\ee
\end{proposition}

\noindent
{\bf Proof.}
This can be proved by using the explicit form of the Killing vectors.

\begin{proposition}
There holds
\bea
\gamma^{AB}\xi_A{}^\alpha\xi_B{}^\mu{}_{;\nu\lambda}
&=& R^\alpha{}_{\lambda\nu\mu}\,.
\label{235}
\eea
\end{proposition}

\noindent
{\bf Proof.}
This follows from eqs. (\ref{22}) and (\ref{267}).

%===========================

\begin{proposition}
There holds
\be
\gamma^{AB}
\xi_A{}^\mu\xi_B{}^\nu{}_{;\beta}=0\,,
\label{2127a}
\ee
\be
\gamma^{AB}
\xi_A{}^\mu{}_{;\alpha}\xi_B{}^\nu{}_{;\beta}
=R^\mu{}_{\alpha}{}^{\nu}{}_{\beta}\,.
\label{242}
\ee
\end{proposition}

\noindent
\noindent\noindent\noindent{\bf Proof.}
Let
\be
\tau^{\mu\alpha}{}_{\nu}
=\gamma^{AB}\xi_A{}^\mu\xi_B{}^\alpha{}_{;\nu}\,
\ee
and
\be
\theta^{\mu}{}_{\alpha}{}^{\nu}{}_{\beta}
=\gamma^{AB}
\xi_A{}^\mu{}_{;\alpha}\xi_B{}^\nu{}_{;\beta}\,. 
\ee
We compute
\be
\nabla_\beta\tau_{\mu\alpha\nu}
=\theta_{\mu\beta\nu\alpha}
-R_{\mu\beta\nu\alpha}\,,
\label{245}
\ee
and
\be
\nabla_\gamma\nabla_\beta\tau_{\mu\alpha\nu}
=R_{\mu\beta\gamma\rho}\tau^{\rho}{}_{\nu\alpha}
+R_{\nu\alpha\gamma\rho}\tau^{\rho}{}_{\mu\beta}\,.
\ee
All higher derivatives of $\tau_{\mu\nu\alpha}$ are expressed
linearly in terms of $\tau_{\mu\nu\alpha}$ and its first
derivative $\nabla_\beta\tau_{\mu\alpha\nu}$ with coefficients
polynomial in curvature.

Let $x'$ be a fixed point. We will show that the tensor 
$\tau_{\mu\nu\alpha}$ together with all its covariant derivatives is
equal to zero at $x=x'$. This will then mean that
$\tau_{\mu\nu\alpha}=0$ identically and, therefore, from eq. (\ref{245})
that $\theta^\mu{}_{\alpha\nu\beta}=R^\mu{}_{\alpha\nu\beta}$.

We have
\be
\tau^{\mu\alpha}{}_{\nu}
=\delta^{ab}P_a{}^\mu P_b{}^\alpha{}_{;\nu}
+\beta^{ij}L_i{}^\mu L_j{}^\alpha{}_{;\nu}\,.
\ee
and
\be
\theta^\mu{}_{\alpha}{}^{\nu}{}_{\beta}
=\delta^{ab} P_a{}^\mu{}_{;\alpha} P_b{}^\nu{}_{;\beta}
+\beta^{ij} L_i{}^\mu{}_{;\alpha} L_j{}^\nu{}_{;\beta}\,.
\ee
Therefore,
\be
\tau^\mu{}_{\alpha\nu}\big|_{x=x'}=0\,
\ee
and
\be
\theta^\mu{}_{\alpha\nu\beta}\big|_{x=x'}
=R^\mu{}_{\alpha\nu\beta}\,.
\ee
Therefore,
\be
\nabla_\beta\tau^\mu{}_{\alpha\nu}\big|_{x=x'}=0\,.
\ee
Thus, by induction, all derivatives of $\tau_{\mu\nu\alpha}$
vanish, and, therefore, $\tau_{\mu\nu\alpha}=0$
identically.
This also proves (\ref{242}) by making use of (\ref{245}).

%=========================================================

Let $i,j$ be non-negative integers.
We define the tensors $X_{(i,j)}$ which are bilinear in Killing 
vectors by
\be
X_{(i,j)}{}^{\mu\nu}{}_{\alpha_1\dots\alpha_i\beta_1\dots\beta_j}
=\gamma^{AB}\nabla_{\alpha_1}\cdots\nabla_{\alpha_i}\xi_A{}^\mu
\nabla_{\beta_1}\cdots\nabla_{\beta_j}\xi_B{}^\nu\,.
\ee

\begin{theorem}
\begin{enumerate}
\item
The tensors 
$X_{(i,j)}$ are $G$-invariant and parallel, that is,
\be
\nabla_\lambda 
X_{(i,j)}
=0\,.
\ee
\item
For even $(i+j)$ the tensors 
$X_{(i,j)}$
are polynomial in the curvature tensor.
\item
For odd $(i+j)$ the tensors 
$X_{(i,j)}$
are identically equal to zero.
\end{enumerate}
\end{theorem}

\noindent
\noindent\noindent\noindent{\bf Proof.} First of all, we notice that (1)
follows from (2) and (3).

There are three cases: a) both $i=2k$ and $j=2m$ are even, b) both
$i=2k+1$ and $j=2m+1$ are odd, and c) $i=2k$ is even and $j=2m+1$ is
odd.

In the case (a), when both $i$ and $j$ are even, by using the  eqs.
(\ref{250}) and (\ref{267}) we immediately obtain a polynomial in the
curvature.

In the cases (b) and (c)  by using the eqs. (\ref{250}) and (\ref{251})
we reduce it to the tensors
$\gamma^{AB}\xi_A{}^\mu{}_{;\alpha}\xi_B{}^\nu{}_{;\beta}$ and
$\gamma^{AB}\xi_A{}^\mu\xi_B{}^\nu{}_{;\beta}\,$. Now, by using the
lemma we prove the theorem.

\begin{proposition}
There holds
\be
\gamma^{AB}\xi_A{}^\mu\eta_B{}^i=0\,,
\label{2140a}
\ee
\be
\gamma^{AB}\eta_A{}^i\eta_B{}^j=\beta^{ij}\,.
\label{2141a}
\ee
\end{proposition}

\noindent\noindent{\bf Proof.} 
This follows from the definition of $\eta_A{}^i$
(\ref{2116}) and eqs. (\ref{2127a}) and (\ref{242}).

%===============================================================
%==============================================================
\subsection{Homogeneous Vector Bundles}

Equation (\ref{227}) imposes strong constraints on the curvature of the
homogeneous bundle ${\cal W}$. We define
\be
{\cal B}_{ab}={\cal F}_{cd}q^c{}_a q^d{}_b\,,
\label{2139}
\ee
\be
{\cal E}_{ab}={\cal F}_{cd}h^c{}_bh^d{}_b\,,
\ee
so that
\be
{\cal B}_{ab}h^a{}_c=0\,,\qquad
{\cal E}_{ab}q^a{}_c=0\,.
\ee

Then, from eq. (\ref{227}) we obtain
\be
[{\cal B}_{ab}, {\cal B}_{cd}]=[{\cal B}_{ab}, {\cal E}_{cd}]=0\,,
\ee
and
\be
[{\cal E}_{cd}, {\cal E}_{ab}] 
- R^f{}_{acd}{\cal E}_{fb}
- R^f{}_{bcd}{\cal E}_{af}
= 0\,.
\label{227a}
\ee
This means that  ${\cal B}_{ab}$ takes values in an Abelian ideal of the
gauge algebra ${\cal G}_{YM}$ and ${\cal E}_{ab}$ takes values in the
holonomy algebra. More precisely, eq. (\ref{227a}) is only possible if
the holonomy algebra ${\cal H}$ is an ideal of the gauge algebra ${\cal
G}_{YM}$. Thus, the gauge group $G_{YM}$ must have a subgroup
$Z\times H$, where $Z$ is an Abelian group and $H$ is the holonomy
group.

We proceed in the following way. The matrices $D^a{}_{ib}$ provide a
natural embedding of the holonomy algebra ${\cal H}$ in the orthogonal
algebra ${\cal SO}(n)$ in the following sense. Let $X_{ab}$ be the
generators of the orthogonal algebra ${\cal SO}(n)$ is some
representation satisfying the commutation relations (\ref{280}).
Let $T_i$ be the matrices defined by
\be
T_i=-\frac{1}{2}D^a{}_{ib}X^b{}_a\,.
\ee

\begin{proposition}
The matrices $T_i$ satisfy the commutation relations
\be
[T_i, T_k]=F^j{}_{ik}T_j\,
\label{2127}
\ee
and form a representation $T$ of the holonomy algebra
${\cal H}$.
\end{proposition}
This can be proved by taking into account  the orthogonal algebra
(\ref{280}).

Thus $T_i$ are the generators of the gauge algebra ${\cal G}_{YM}$
realizing a representation $T$ of the holonomy algebra ${\cal H}$. 
Since ${\cal B}_{ab}$ takes values in the Abelian ideal of the algebra
of the gauge group we also have
\be
[{\cal B}_{ab},T_j]=0\,.
\ee

Then by using eq. (\ref{313}) one can show that
\footnote{We correct here a
sign misprint in eq. (3.24) in \cite{avramidi96}.}
\be
{\cal E}_{ab}=\frac{1}{2}R^{cd}{}_{ab}X_{cd}
=-E^i{}_{ab}T_i\,.
\ee

\begin{proposition}
The two form
\bea
{\cal F}_{ab}&=&
-E^i{}_{ab}T_i+{\cal B}_{ab}
\nonumber\\
&=&\frac{1}{2}R^{cd}{}_{ab}X_{cd}+{\cal B}_{ab}
\,
\label{2153mm}
\eea
satisfies the constrains (\ref{227}), and, therefore, gives the
curvature of the homogeneous bundle ${\cal W}$.
\end{proposition}

Now, we consider the representation $\Sigma$ of the orthogonal algebra
defining the spin-tensor bundle ${\cal T}$ and define the matrices
\be
G_{ab}=\Sigma_{ab}\otimes\II_X+\II_\Sigma\otimes X_{ab}\,.
\ee
Obviously, these matrices are the generators of the orthogonal algebra
in the product representation $\Sigma\otimes X$.

Next, the matrices 
\be
Q_i=-\frac{1}{2}
D^a{}_{ib}\Sigma^b{}_a\,
\label{2104}
\ee
form a representation $Q$ of the holonomy algebra
${\cal H}$.
and the matrices 
\bea
{\cal R}_{i}&=&
Q_i\otimes\II_{T}+\II_{\Sigma}\otimes T_i
\nonumber\\
&=&
-\frac{1}{2}D^a{}_{ib}G^b{}_a
\,
\label{2156mm}
\eea
are the generators of the holonomy algebra in the
product representation ${\cal R}=Q\otimes T$. 

Then the total curvature,
that is, the commutator of covariant derivatives, (\ref{220}) of  a
twisted spin-tensor bundle ${\cal V}$ is
\bea
\Omega_{ab}&=&-E^i{}_{ab}{\cal R}_i
+{\cal B}_{ab}
\nonumber\\
&=&\frac{1}{2}R^{cd}{}_{ab}G_{cd}+{\cal B}_{ab}
\,.
\label{2144}
\eea

Finally, we define the Casimir operators of the holonomy algebra in
the representations $Q$, $T$ and ${\cal R}$
\be
T^2=C_2(H,T)=\beta^{ij}T_i T_j
=\frac{1}{4}R^{abcd}X_{ab}X_{cd}\,,
\ee
\be
%\qquad
Q^2=C_2(H,Q)=\beta^{ij}Q_i Q_j
=\frac{1}{4}R^{abcd}\Sigma_{ab}\Sigma_{cd}\,,
\ee
\be
{\cal R}^2=C_2(H,{\cal R})=\beta^{ij}{\cal R}_i {\cal R}_j
=\frac{1}{4}R^{abcd}G_{ab}G_{cd}\,.
\ee
They commute with all matrices $T_i$, $Q_i$ and ${\cal R}_i$
respectively.

%===========================================================
\subsection{Twisted Lie Derivatives}

Let $\varphi$ be a section of a  twisted homogeneous spin-tensor bundle
${\cal T}$. Let $\xi_A$ be the basis of Killing vector fields. Then the
covariant (or generalized, or twisted)  
Lie derivative of $\varphi$ along $\xi_A$ is defined by
\be
{\cal L}_A\varphi=
{\cal L}_{\xi_A}\varphi
=\left(\nabla_{\xi_A}+S_A\right)\varphi\,,
\label{2160x}
\ee
where $\nabla_{\xi_A}=\xi_A{}^\mu\nabla_\mu$,
\be
S_A=\eta_A{}^i{\cal R}_i
=\frac{1}{2}\xi_A{}^a{}_{;b}G^b{}_a\,,
\ee
and $\eta_A{}^i$ are defined by (\ref{2116}).
Note that
\be
S_a q^a{}_b=0\,.
\ee

\begin{proposition}
There hold
\be
[\nabla_{\xi_A},\nabla_{\xi_B}]\varphi
=\left(
C^C{}_{AB}\nabla_{\xi_C}-{\cal R}_{AB}
+{\cal B}_{AB}
\right)\varphi\,,
\label{2160}
\ee
\be
\nabla_{\xi_A}S_B
={\cal R}_{AB}\,,
\label{2152}
\ee
\be
[S_A,S_B]=C^C{}_{AB}S_C
-{\cal R}_{AB}\,,
\label{2165}
\ee
where
\bea
{\cal R}_{AB}&=&
\xi_A{}^a\xi_B{}^bE^i{}_{ab} {\cal R}_i
\nonumber\\
&=&
-\frac{1}{2}R^{cd}{}_{ab}\xi_A{}^a\xi_B{}^b G_{cd}
\,,
\eea
\be
{\cal B}_{AB}=\xi_A{}^a\xi_B{}^b{\cal B}_{ab}\,.
\ee
\end{proposition}

\noindent
{\bf Proof.}
By using the properties of the Killing vectors described in the previous
section and the eq. (\ref{2144})
we obtain first (\ref{2160}).
Next, by using the eqs. (\ref{2117}) we obtain (\ref{2152}), and,
further, by using the eq. (\ref{2124}) we get (\ref{2165}).

Notice that from the definition (\ref{2139}) we have
\be
P_c{}^a L_i{}^b{\cal B}_{ab}=
L_i{}^a L_j{}^b{\cal B}_{ab}=0\,,
\label{2169}
\ee
and
\be
P_c{}^a P_d{}^b{\cal B}_{ab}={\cal B}_{cd}\,.
\ee
This means that the matrix ${\cal B}_{AB}$ has the form
\be
{\cal B}_{AB}=\left(
\begin{array}{cc}
{\cal B}_{ab} & 0\\
0 & 0\\
\end{array}
\right)\,,
\label{2171xx}
\ee
and, therefore,
\be
C^A{}_{BC}{\cal B}_{AD}=\gamma^{BD}C^A{}_{BC}{\cal B}_{DE}=0\,.
\label{2172}
\ee

We define the  operator
\be
{\cal L}^2=\gamma^{AB}{\cal L}_A{\cal L}_B\,.
\ee

\begin{theorem}
The operators ${\cal L}_A$ and ${\cal L}^2$
satisfy the commutation relations
\be
[{\cal L}_A, {\cal L}_B] = C^C{}_{AB} {\cal L}_C+{\cal B}_{AB},
\label{2173xx}
\ee
or, in more details,
\bea
[{\cal L}_a, {\cal L}_b] &=& E^i{}_{a b} {\cal L}_i+{\cal B}_{ab},
\\{}
[{\cal L}_i,{\cal L}_a] &=& D^b{}_{ia}{\cal L}_b,
\\{}
[{\cal L}_i, {\cal L}_j] &=& F^k{}_{ij} {\cal L}_k\,,
\label{321b}
\eea
and
\be
[{\cal L}_A,{\cal L}^2]=2\gamma^{BC}{\cal B}_{AB}{\cal L}_C\,.
\label{2178xx}
\ee
\end{theorem}

\noindent
\noindent\noindent\noindent{\bf Proof.} 
This follows from 
\be
[{\cal L}_A,{\cal L}_B]=[\nabla_{\xi_A},\nabla_{\xi_B}]
+[\nabla_{\xi_A},S_B]
-[\nabla_{\xi_B},S_A]
+[S_A,S_B]
\ee
and eqs. (\ref{2160}), (\ref{2152}), and (\ref{2165}). The eq.
(\ref{2178xx}) follows directly from (\ref{2173xx}).

The operators ${\cal L}_A$ form an algebra that
is a direct sum of a nilpotent ideal and a 
semisimple algebra. For the lack of a better name we call this algebra
{\it gauged curvature algebra} and denote it by ${\cal G}_{\rm gauge}$.

\begin{proposition}
There hold
\be
\gamma^{AB}\xi_A{}^\mu S_B=0\,,
%\qquad
\ee
\be
\gamma^{AB}\nabla_{\xi_A} S_{B}=0\,,
\label{2153}
\ee
\be
\gamma^{AB}S_A S_B={\cal R}^2\,.
\label{2154}
\ee
\end{proposition}

\noindent
{\bf Proof.} This can be proved by using the eqs. (\ref{2140a}),
(\ref{2152}) and (\ref{2141a}).

\begin{theorem}
The Laplacian $\Delta$ acting on sections of a twisted spin-tensor
bundle ${\cal V}$ over a symmetric space has the form
\be
\Delta={\cal L}^2-{\cal R}^2\,.
\label{2110}
\ee

\be
[{\cal L}_{A},\Delta] = 2 \gamma^{BC}{\cal B}_{AB}{\cal L}_C\,.
\label{43}
\ee
\end{theorem}

\noindent
\noindent\noindent\noindent{\bf Proof.} We have
\be
\gamma^{AB}{\cal L}_A{\cal L}_B=\gamma^{AB}\nabla_{\xi_A}\nabla_{\xi_B}
+\gamma^{AB}S_A\nabla_{\xi_B}
+\gamma^{AB}\nabla_{\xi_A}S_B
+\gamma^{AB}S_AS_B\,.
\ee
Now, by using eqs. (\ref{267}) and (\ref{2127a}) we get
\be
\gamma^{AB}\nabla_{\xi_A}\nabla_{\xi_B}=\Delta\,.
\ee
Next, by using the eqs. (\ref{2152}), (\ref{2153}) and (\ref{2154}), we
obtain (\ref{2110}). The eq. (\ref{43}) follows from the commutation
relations  (\ref{2173xx}).

%======================================================
\subsection{Isometries and Pullbacks}

Let $\omega^i$ be the canonical coordinates on the holonomy group and 
$(k^A)=(p^a,\omega^i)$ be the canonical coordinates on the  gauged
curvature group.  We fix a point $x'$ so that the basis Killing vectors
fields $\xi_A$ satisfy the initial  conditions (\ref{237})-(\ref{238}) 
and are given by (\ref{226})-(\ref{229}). Let
$\xi=\left<k,\xi\right>=k^A\xi_A
=p^a P_a+\omega^iL_i$ be a Killing vector field
and let $\psi_t: M\to M$ be the
one-parameter diffeomorphism (the isometry) generated by the vector
field  $\xi$. Let $\hat x=\psi_t(x)$, so that
\be
\frac{d\hat x}{dt}=\xi{}(\hat x)\,
\label{610}
\ee
and
\be
\hat x\big|_{t=0}=x\,.
\label{2189xx}
\ee
The solution of this equation depends on the parameters 
$t,p,\omega,x$ and $x'$, that is, 
\be
\hat x=\hat x(t,p,\omega,x,x')\,.
\label{2190nnn}
\ee

We will be interested mainly in the case when the points $x$ and $x'$
are close to each other. In fact, at the end of our calculations
we will take the limit $x=x'$.
In this case, as we will show below, the Jacobian 
\be
\det\left(\frac{\partial \hat x^\mu}{\partial p^a}\right)\ne 0
\ee
is not equal to zero,
and, therefore, coordinates $p$ can be used to parametrize the 
point $\hat x$, that is, the eq. (\ref{2190nnn}) defines the function
\be
p=p(t,\omega,\hat x, x,x')\,.
\ee
We will be interested in those trajectories that reach the point 
$x'$ at the time $t=1$.
So, we look at the values $\hat x(1,p,\omega,x,x')$
when the parameters $p$ are varied. Then, as we will show below,
there is always a value of the parameters $p$ that we call $\bar p$
such that 
\be
\hat x(1,\bar p,\omega,x,x')=x'\,.
\label{2190nn}
\ee
Thus, eq. (\ref{2190nn}) defines a function $\bar p=\bar p(\omega,x,x')$.
Therefore, the parameters $\bar p$ can be used to parameterize
the point $x$. Of course,
\be
\bar p(\omega,x,x')=p(1,\omega,x',x,x')\,.
\ee

Now, we choose the normal coordinates $y^a$ of the point   defined above
and the normal coordinates $\hat y^a$ of the point $\hat x$  with the
origin at $x'$, so that  the normal coordinates $y'$ of the point $x'$
are  equal to zero, $y'^a=0$. Recall that the normal coordinates are
equal to the  components of the tangent vector at the point $x'$  to the
geodesic connecting the points $x'$ and the current point, that is,
$y^a=-e^a{}_{\mu'}(x')\sigma^{;\mu'}(x,x')$ and $\hat
y^a=-e^a{}_{\mu'}(x')\sigma^{;\mu'}(\hat x,x')$. Then by taking into
account eqs. (\ref{226}) and (\ref{229}) the equation (\ref{610})
becomes
\be
{d \hat y^a\over dt} = \left(\sqrt {K(\hat y)}\cot \sqrt
{K(\hat y)}\right)^a_{\ b}p^b -\omega^iD^a{}_{ib}\hat y^b\,,
\label{447}
\ee
with the initial condition
\be
\hat y^a\big|_{t=0}=y^a\,.
\label{2189}
\ee
The solution of this equation defines a function 
$\hat y=\hat y(t,p,\omega,y)$. 

\begin{proposition}
The Taylor expansion of the solution 
$\hat y=\hat y(t,p,\omega,y)$ of the eq. (\ref{447})
in $t$ reads
\be
\hat y^a=y^a
+\left[\left(\sqrt {K(y)}\cot \sqrt
{K(y)}\right)^a_{\ b}p^b -\omega^iD^a{}_{ib}y^b\right]t
+O(t^2)\,.
\label{447nn}
\ee
The Taylor expansion of the function $\hat y=\hat y(t,p,\omega,y)$
in $p$ and $y$ reads
\be
\hat y^a = \left(\exp[-tD(\omega)]\right)^a{}_b y^b
+\left({1-\exp[-t D(\omega)]\over D(\omega)}\right)^a{}_bp^b
+O(y^2,p^2,py)\,.
\label{2197mm}
\ee
There holds
\be
\det\left(
\frac{\partial \hat y^a}{\partial p^b}
\right)\Bigg|_{p=y=0,t=1}
=\det{}_{TM}\left({\sinh[\,D(\omega)/2]\over D(\omega)/2}\right)\,.
\label{429a}
\ee
\end{proposition}
{\bf Proof.}
The eq. (\ref{447nn}) trivially follows from the eq. 
(\ref{447}).

Let us expand the function $\hat y(t,p,\omega,y)$
in Taylor series in $p$ and $y$
restricting ourselves to linear terms,
that is,
\be
\hat y^a = \hat y^a\Big|_{p=y=0}
+\frac{\partial \hat y^a}{\partial p^b}\Big|_{p=y=0}p^b
+\frac{\partial \hat y^a}{\partial y^b}\Big|_{p=y=0}y^b
+O(p^2,y^2,py)\,.
\label{2199mm}
\ee
First of all, for $p=0$ the eq. (\ref{447}) becomes
\be
{d \hat y^a\over dt} = -\omega^iD^a{}_{ib}\hat y^b\,.
\label{447mm}
\ee
The solution of this equation with the initial condition
$\hat y=0$ is trivial, therefore,
\be
\hat y\Big|_{p=y=0}=\hat y(t,0,\omega,0)=0\,.
\label{2203mm}
\ee
Next, by differentiating the eq. (\ref{447mm}) with respect to 
$y^b$ and setting $y=0$ we obtain the equation
\be
\frac{d}{dt}\frac{\partial \hat y^a}{\partial y^b}\Big|_{p=y=0}
 = -\omega^iD^a{}_{ic}
\frac{\partial \hat y^c}{\partial y^b}\Big|_{p=y=0}\,.
\label{447cc}
\ee
with the initial condition
\be
\frac{\partial \hat y^a}{\partial y^b}\Big|_{p=y=t=0}
=\delta^a{}_b\,.
\ee
The solution of this equation is 
\be
\frac{\partial \hat y^a}{\partial y^b}\Bigg|_{p=y=0}
=\left(\exp[-tD(\omega)]\right)^a{}_b\,,
\label{2204mm}
\ee
where
\be
D(\omega)=\omega^i D_i\,.
\ee

Let
\be
Z^a{}_b=\frac{\partial \hat y^a}{\partial p^b}\Big|_{p=y=0}\,.
\ee
Then by differentiating the eq. (\ref{447}) with respect to $p^b$
and setting $p=0$, we obtain
\be
{d \over dt}Z^a{}_b = \delta^a{}_b
-\omega^iD^a{}_{ic}Z^c{}_b\,,
\label{447b}
\ee
with the initial condition
\be
Z^a{}_b\big|_{t=0}=0\,.
\ee
The solution of this equation is
\be
Z={1-\exp[-t D(\omega)]\over D(\omega)}\,.
\label{452}
\ee
By substituting the eqs. (\ref{2203mm}), (\ref{2204mm})
and (\ref{452}) in (\ref{2199mm}) we get the desired
result (\ref{2197mm}).

Finally, 
by taking into account that the matrix $D(\omega)$ is traceless,
we find first 
$\det\exp[tD(\omega)]=1$, and, then
by using eq. (\ref{452}) we obtain (\ref{429a}).

%============
The function $\hat y=\hat y(t,p,\omega,y)$
implicitly defines the function
\be
p=p(t,\omega,\hat y,y)\,.
\ee
The function 
$\bar p=\bar p(\omega,y)$ is now defined
by the equation
\be
\hat y(1,\bar p,\omega,y)=0\,,
\label{2197nn}
\ee
or
\be
\bar p(\omega,y)=p(1,\omega,0,y)\,.
\ee

\begin{proposition}
The Taylor expansion of the function $\bar p(\omega,y)$
in $y$ has the form
\be
\bar p^a=-\left(D(\omega)
\frac{\exp[-D(\omega)]}{1-\exp[-D(\omega)]}
\right)^a{}_b y^b
+O(y^2)\,.
\label{2213mm}
\ee
Therefore,
\be
\det\left(-
\frac{\partial \bar p^a}{\partial y^b}
\right)\Bigg|_{y=0}
=\det{}_{TM}\left({\sinh[\,D(\omega)/2]
\over D(\omega)/2}\right)^{-1}\,.
\label{429aa}
\ee
\end{proposition}
{\bf Proof.}
We expand $\bar p$ in Taylor series in $y$
\be
\bar p^a = \bar p^a\Big|_{y=0}
+\frac{\partial \bar p^a}{\partial y^b}\Big|_{y=0}y^b
+O(y^2)\,.
\ee
Next, by taking into account (\ref{2203mm}) we have
\be
\bar p\Big|_{y=0}=0\,.
\label{2216mm}
\ee
Further, by differentiating (\ref{2197nn}) with respect to
$y^c$ and setting $y=0$ we get
\be
\frac{\partial \hat y^a}{\partial y^b}\Big|_{p=y=0, t=1}
+\frac{\partial \hat y^a}{\partial p^c}\Big|_{p=y=0, t=1}
\frac{\partial \bar p^c}{\partial y^b}\Big|_{y=0}=0\,,
\ee
and, therefore,
\be
\frac{\partial \bar p^a}{\partial y^b}\Big|_{y=0}
=-\left(D(\omega)
\frac{\exp[-D(\omega)]}{1-\exp[-D(\omega)]}
\right)^a{}_b\,.
\ee
This leads to both (\ref{2213mm}) and (\ref{429aa}).

%===================================
Now, we define
\be
\Lambda^{\hat\mu}{}_{\nu}=
\frac{\partial \hat x^\mu}{\partial x^{\nu}}\,.
\ee
The pullback of the metric by the diffeomorphism $\psi_t$ is defined by
\be
(\psi_t^* g)_{\mu\nu}(x)=
\Lambda^{\hat\alpha}{}_{\mu}\Lambda^{\hat\beta}{}_{\nu}
g_{\hat\alpha\hat\beta}(\hat x)\,.
\ee
Since $\psi_t$ is an isometry, we have
\be
(\psi_t^* g)_{\mu\nu}(x)
=g_{\mu\nu}(x)\,.
\label{611}
\ee
Therefore, the inverse matrix $\Lambda^{-1}$ is equal to
\be
(\Lambda^{-1})^{\mu}{}_{\hat\alpha}
=g^{\mu\nu}(x)\Lambda^{\hat\beta}{}_{\nu}
g_{\hat\beta\hat\alpha}(\hat x)\,.
\ee
Let $e^a{}_\mu$ and $e_a{}^\mu$ be a local orthonormal frame 
that is obtained by parallel transport along geodesics from a
point $x'$.
Then 
the action of the pullback $\psi_t^*$ on the orthonormal frame
is
\be
(\psi_t^*e^a{})_{\mu}(x)=\Lambda^{\hat\alpha}{}_{\mu} 
e^a{}_{\hat\alpha}(\hat x)\,.
\ee
Since $\psi_t$ is an isometry, we have
\be
\delta_{ab}(\psi_t^* e^a)_{\alpha}(x)
(\psi_t^* e^b)_{\beta}(x)
=\delta_{ab}e^a{}_{\alpha}(x) e^b{}_{\beta}(x)\,.
\ee
Therefore, the frames 
of $1$-forms $e^a$ and $\psi_t^* e^a$ are related by an 
orthogonal transformation
\be
(\psi_t^* e^a)(x)=O^a{}_b e^b(x)\,,
\ee
where the matrix $O^a{}_b$ is defined by
\be
O^a{}_b
=e^a{}_{\hat\alpha}(\hat x)
\Lambda^{\hat\alpha}{}_{\mu} e_b{}^{\mu}(x)
\,.
\ee

\begin{proposition}
For $p=y=0$ the matrix $O$ has the form
\be
O\Big|_{p=y=0}=\exp\left[-t D(\omega)\right]\,.
\label{2149}
\ee
\end{proposition}
\noindent\noindent\noindent{\bf Proof.}
We use normal coordinates $\hat y^a$ and $y^a$. Then the matrix
$O$ takes the form
\be
O^a{}_b=e^a{}_{\hat\alpha}
\frac{\partial \hat x^\alpha}{\partial \hat y^c}
\frac{\partial \hat y^c}{\partial y^d}
\frac{\partial y^d}{\partial x^\mu}e_b{}^\mu\,.
\ee
Now, by using the Jacobian matrix (\ref{244mm}) and recalling
that $\hat y=0$ for $p=y=0$ we obtain
\be
\frac{\partial y^a}{\partial x^\mu}e_b{}^\mu
\Bigg|_{p=y=0}
=e^a{}_{\hat\alpha}
\frac{\partial \hat x^\alpha}{\partial \hat y^b}
\Bigg|_{p=y=0}
=\delta^a{}_b\,.
\ee
Therefore,
\be
O^a{}_b\Big|_{p=y=0}
=\frac{\partial \hat y^a}{\partial y^b}\Bigg|_{p=y=0}\,,
\ee
and, finally (\ref{2204mm}) gives the desired result
(\ref{2149}).

Let $\varphi$ be a section of the twisted spin-tensor bundle ${\cal V}$.
Let $V_{x}$ be the fiber at the point $x$ and
$V_{\hat x}$ be the fiber at the point $\hat x=\psi_t(x)$.
The pullback of the diffeomorphism $\psi_t$ defines the map,
that we call just the pullback,
\be
\psi_t^*: C^\infty({\cal V})\to C^\infty(\cal V)
\ee
on smooth sections of the twisted spin-tensor bundle ${\cal V}$.
The pullback of tensor fields of type $(p,q)$ is
defined by
\be
(\psi_t^*\varphi)^{\mu_1\dots \mu_p}_{\nu_1\dots\nu_q}(x)
=\Lambda^{\hat\beta_1}{}_{\nu_1}\cdots
\Lambda^{\hat\beta_q}{}_{\nu_q}
(\Lambda^{-1})^{\mu_1}{}_{\hat\alpha_1}\cdots
(\Lambda^{-1})^{\mu_p}{}_{\hat\alpha_p}
\varphi^{\hat\alpha_1\dots\hat\alpha_p}_{\hat\beta_1\dots\hat\beta_q}(\hat x)\,.
\ee
We define the twisted pullback (a combination of a proper pullback and a
gauge transformation)  of a tensor of type $(p,q)$ by
\be
(\psi_t^*\varphi)^{a_1\dots a_p}_{b_1\dots b_q}(x)
=O^{c_1}{}_{b_1}\cdots O^{c_q}{}_{b_q}
O_{d_1}{}^{a_1}\cdots O_{d_p}{}^{a_p}
\varphi^{d_1\dots d_p}_{c_1\dots c_q}(\hat x)\,.
\ee
Since the matrix $O$ is orthogonal, it can be
parametrized by
\be
O=\exp\theta\,,
\label{2234mm}
\ee
where $\theta_{ab}$ is an antisymmetric matrix. 
The orthogonal transformation of the frame pulled back 
causes the transformation of spinors
\be
(\psi_t^*\varphi)(x)
=\exp\left(-\frac{1}{4}\theta_{ab}\gamma^{ab}\right)
\varphi(\hat x)\,.
\ee
More generally, we have

\begin{proposition}
Let $\varphi$ be a section of a twisted spin-tensor bundle
{\cal V}. Then
\be
(\psi_t^*\varphi)(x)
=\exp\left(-\frac{1}{2}\theta_{ab}G^{ab}\right)
\varphi(\hat x)\,.
\label{2236mm}
\ee
In particular, for $p=y=0$ (or $x=x'$)
\be
(\psi_t^*\varphi)(x)\Big|_{p=y=0}
=\exp\left[t{\cal R}(\omega)\right]
\varphi(x')\,,
\label{2237mm}
\ee
where
\be
{\cal R}(\omega)=\omega^i{\cal R}_i\,.
\ee
\end{proposition}
{\bf Proof.}
First, from the eq. (\ref{2149}) we see that
\be
\theta^a{}_b\Big|_{p=y=0}=-t\omega^i D^a{}_{ib}\,.
\ee
Then, from the definition (\ref{2156mm})
of the matrices ${\cal R}_i$ we get (\ref{2237mm}).

It is not very difficult to check that
the Lie derivatives are nothing but 
the generators of the pullback, that is,
\be
{\cal L}_{\xi}\varphi
=k^A{\cal L}_A\varphi
=\frac{d}{dt}(\psi_t^*\varphi)
\Big|_{t=0}\,.
\ee
We will use this fundamental 
fact to compute the heat kernel diagonal
below.

%===============================================================
\section{Heat Semigroup}
\setcounter{equation}0

%============================================
\subsection{Geometry of the Curvature Group}

Let $G_{\rm gauge}$ be the gauged curvature group and $H$ be its
holonomy subgroup.  Both these groups have compact algebras. However,
while the holonomy group is always compact, the curvature group is, in
general, a product of a nilpotent group, $G_0$, and a
semi-simple group, $G_s$,
\be
G_{\rm gauge}=G_0\times G_s\,.
\ee
The semi-simple group $G_s$ is a product $G_s=G_+\times G_-$ of a
compact $G_+$ and a non-compact $G_-$ subgroups.

Let $\xi_A$ be the basis Killing vectors, $k^A$ be the canonical
coordinates on the curvature group $G$ and $\xi(k)=k^A\xi_A$.
The canonical coordinates are exactly the normal coordinates
on the group defined above.
Let $C_A$ be the generators of the curvature group in
adjoint representation and $C(k)=k^AC_A$. 
In the following $\partial_M$ means the partial derivative 
$\partial/\partial k^M$ with
respect to the canonical coordinates.
We define the matrix $Y^A{}_M$ by the equation
\be
\exp[-\xi(k)]\partial_M\exp[\xi(k)]=Y^A{}_M\xi_A\,,
\label{31xx}
\ee
which is well defined since the right hand side lies in the Lie
algebra of the curvature group.
This can be written in the form
\be
\exp[-\xi(k)]\partial_M\exp[\xi(k)]=\exp[-Ad_{\xi(k)}]\partial_M
\label{32xx}
\ee
where the operator $Ad_X$ is defined by
$Ad_X Z=[X,Z]$.
This enables us to compute the matrix $Y=(Y^A{}_M)$ 
explicitly, namely,
\be
Y=\frac{1-\exp[-C(k)]}{C(k)}\,.
\ee
Let $X=(X_A{}^M)=Y^{-1}$ be the inverse matrix of $Y$.
Then we define the $1$-forms $Y^A$ and the vector fields $X_A$
on the group $G$ by
\be
Y^A=Y^A{}_M dk^M\,,\qquad
X_A=X_A{}^M\partial_M\,.
\ee

\begin{proposition}
There holds
\be
X_A \exp[\xi(k)]=\exp[\xi(k)]\xi_A\,.
\label{35xx}
\ee
\end{proposition}
{\bf Proof.} This follows immediately from the eq. (\ref{31xx}).

Next, 
by differentiating the eq. (\ref{31xx}) 
with respect to $k^L$ and alternating the indices $L$ and $M$ we obtain
\be
\partial_L Y^A{}_M-\partial_M Y^A{}_L=-C^A{}_{BC}Y^B{}_{L}Y^C{}_M\,,
\label{36xx}
\ee
which, of course, can also be written as
\be
dY^A=-\frac{1}{2}C^A{}_{BC}Y^B\wedge Y^C\,.
\ee

\begin{proposition}
The vector fields $X_A$ satisfy the commutation relations
\be
[X_A,X_B]=C^C{}_{AB}X_C\,.
\label{37xx}
\ee
\end{proposition}
{\bf Proof.} This follows from the eq. (\ref{36xx}).

The vector fields $X_A$ are nothing but the right-invariant
vector fields. They form a representation  of
the curvature algebra.

We will also need the following fundamental property of Lie groups.

\begin{proposition}
Let $G$ be a Lie group with the structure constants $C^A{}_{BC}$,
$C_A=(C^B{}_{AC})$ and $C(k)=C_Ak^A$. Let $\gamma=(\gamma_{AB})$ be a
symmetric non-degenerate matrix satisfying the  equation
\be
(C_A)^T=-\gamma C_A \gamma^{-1}\,.
\ee
Let $X=(X_A{}^M)$ be a matrix defined by
\be
X=\frac{C(k)}{1-\exp[-C(k)]}\,.
\ee
Then
\be
(\det X)^{-1/2}\gamma^{AB}X_A{}^M\partial_M X_B{}^N\partial_N
(\det X)^{1/2} = -\frac{1}{24}\gamma^{AB}C^C{}_{AD}C^D{}_{BC}\,.
\label{425}
\ee
\end{proposition}
{\bf Proof.} It is easy to check that this equation holds at $k=0$. Now,
it can be proved by showing that it is a group invariant. For a 
detailed proof for semisimple groups see
\cite{helgason84,camporesi90,dowker71}.

It is worth stressing that this equation holds not only on semisimple
Lie groups but on any group with a compact Lie algebra, that is, when
the structure constants $C^A{}_{BC}$ and the matrix $\gamma_{AB}$, used
to define the metric $G_{MN}$ and the operator $X^2$, satisfy the eq.
(\ref{288xx}). Such algebras can have an Abelian center as in eq.
(\ref{290xx}).

Now, by using the right-invariant 
vector fields we define a metric on 
the curvature group $G$
\be
G_{MN} = \gamma_{AB} Y^{A}{}_{M} Y^{B}{}_{N}\,,
\qquad
G^{MN} = \gamma^{AB} X_A{}^{M} X_B{}^{N}\,.
\label{418}
\ee
This metric is bi-invariant and satisfies, in
particular, the equation
\be
{\cal L}_{X_A}G_{BC}=
X_A{}^M\partial_M G_{BC}
+G_{BM}\partial_C X_A{}^M
+G_{MC}\partial_B X_A{}^M=0\,.
\label{311xxx}
\ee
This equation is proved by using eqs. (\ref{288xx}) and 
(\ref{37xx}).
This means that the vector fields $X_A$ are the Killing vector fields
of the metric $G_{MN}$.
One can easily show that 
this metric defines the following natural affine
connection $\nabla^G$ on the group
\be
\nabla^G_{X_C} X_A=-\frac{1}{2}C^A{}_{BC}X_B{}\,,
\qquad
\nabla^G_{X_C} Y^A{}=\frac{1}{2}C^B{}_{AC}Y^B{}\,,
\label{311xx}
\ee
with the scalar curvature
\be
R_G=-\frac{1}{4}\gamma^{AB}C^C{}_{AD}C^D{}_{BC}\,.
\label{312xx}
\ee

Since the matrix $C(k)$ is traceless we have $\det\exp[C(k)/2]=1$,
and, therefore, the volume element on the group is
\be
|G|^{1/2}=\left(\det G_{MN}\right)^{1/2}
=|\gamma|^{1/2}\det{}_{{\cal G}}
\left({\sinh[C(k)/2]\over C(k)/2}\right)\,,
\label{419}
\ee
where $|\gamma|=\det \gamma_{AB}$.
Notice that this function is precisely the inverse Van Vleck-Morette
determinant (\ref{245xx}) on the group
in normal coordinates.

It is not difficult to see that
\be
k^M Y^A{}_M=k^M X_M{}^A=k^A\,.
\label{314xx}
\ee
By differentiating this equation with respect to $k^B$ and 
contracting the indices $A$ and $B$ we obtain
\be
k^M\partial_A X_M{}^A=N-X_A{}^A\,.
\label{315xx}
\ee

Now, by contracting the eq. (\ref{311xx}) with $G^{BC}$ we obtain
the zero-divergence condition for the right-invariant vector fields
\be
|G|^{-1/2}\partial_M\left(|G|^{1/2}X_A{}^M\right)=0\,.
\label{313xx}
\ee
Next, we define the Casimir operator
\be
X^2=C_2(G,X)=\gamma^{AB}X_AX_B\,.
\ee
By using the eq. (\ref{313xx}) one can easily show that $X^2$ is
an invariant differential operator that is nothing
but the scalar Laplacian on the group
\be
X^2=|G|^{-1/2}\partial_M |G|^{1/2}G^{MN}\partial_N
=G^{MN}\nabla^G_M\nabla^G_N\,.
\ee
Then, by using the eqs. (\ref{288xx}) and (\ref{320}) 
one can show that the operator $X^2$ commutes with the operators
$X_A$,
\be
[X_A,X^2]=0\,.
\ee

Since we will actually be working with the gauged curvature group,
we introduce now the operators (covariant right-invariant
vector fields) $J_A$ by
\be
J_A=X_A-\frac{1}{2}{\cal B}_{AB}k^B\,,
\ee
and the operator 
\be
J^2=\gamma^{AB}J_AJ_B\,.
\ee

\begin{proposition}
The operators $J_A$ and $J^2$ satisfy the commutation relations
\be
[J_A,J_B]=C^C{}_{AB}J_C+{\cal B}_{AB}\,,
\label{326xx}
\ee
and
\be
[J_A,J^2]=2{\cal B}_{AB}J^B\,.
\label{326xxx}
\ee
\end{proposition}
{\bf Proof.}
By using the eqs. 
(\ref{2169})-(\ref{2172})
we obtain
\be
X_B{}^A {\cal B}_{AM}
=
\gamma_{BN}\gamma^{AC}X_C{}^N {\cal B}_{AM}
={\cal B}_{BM}\,,
\label{327xx}
\ee
and, hence,
\be
\gamma^{AB}X_B{}^M {\cal B}_{AM}
=0\,,
\ee
and, further, by using (\ref{37xx})
we obtain (\ref{326xx}).
By using the eqs. (\ref{327xx}) we get (\ref{326xxx}).

Thus, the operators $J_A$ form a representation of the gauged curvature
algebra.
Now, let ${\cal L}_A$ be the operators of Lie derivatives
satisfying the commutation relations (\ref{2173xx}) and
${\cal L}(k)=k^A {\cal L}_A$. 

\begin{proposition}
There holds
\be
J_A \exp[{\cal L}(k)]=\exp[{\cal L}(k)]{\cal L}_A\,.
\ee
and, therefore,
\be
J^2 \exp[{\cal L}(k)]=\exp[{\cal L}(k)]{\cal L}^2\,.
\label{330xx}
\ee
\end{proposition}
{\bf Proof.}
Similarly to (\ref{32xx})
we have
\be
\exp[-{\cal L}(k)]\partial_M\exp[{\cal L}(k)]
=\exp[-Ad_{{\cal L}(k)}]\partial_M\,.
\ee
By using the commutation relations (\ref{2173xx}) and eq. (\ref{2172})
we obtain
\be
\exp[-{\cal L}(k)]\partial_M\exp[{\cal L}(k)]
=Y^A{}_M{\cal L}_A+\frac{1}{2}{\cal B}_{MN}k^N\,.
\ee
The statement of the proposition follows from the definition of the
operators $J_A$, $J^2$ and ${\cal L}^2$.

%==============================================
\subsection{Heat Kernel on the Curvature Group}

Let ${\cal B}$ be the matrix with the components ${\cal
B}=(\gamma^{AB}{\cal B}_{BC})$.
Let $k^A$ be the canonical coordinates on the curvature group $G$
and $A(t;k)$ be a function defined by
\be
A(t;k)=\det{}_{\cal G}
\left(\frac{\sinh\left[C(k)/2+t{\cal B}\right]}
{C(k)/2+t{\cal B}}\right)^{-1/2}\,.
\ee
By using the eqs. (\ref{327xx}) one can rewrite this in the form
\be
A(t;k)=\det{}_{\cal G}
\left(\frac{\sinh\left[C(k)/2\right]}
{C(k)/2}\right)^{-1/2}
\det{}_{\cal G}
\left(\frac{\sinh\left[t{\cal B}\right]}
{t{\cal B}}\right)^{-1/2}
\,.
\label{334xx}
\ee
Notice also that due to (\ref{2171xx})
\be
\det{}_{\cal G}
\left(\frac{\sinh\left[t{\cal B}\right]}
{t{\cal B}}\right)^{-1/2}
=\det{}_{TM}
\left(\frac{\sinh\left[t{\cal B}\right]}
{t{\cal B}}\right)^{-1/2}\,,
\ee
where ${\cal B}$ is now regarded as just the matrix 
${\cal B}=({\cal B}^a{}_b)$.

Let $\Theta(t;k)$ be another function on the group $G$ defined by
\be
\Theta(t;k)=\frac{1}{2}\left<k, \gamma\hat\Theta k\right>\,,
\label{337mm}
\ee
where $\hat\Theta$ is the matrix
\be
\hat\Theta=t{\cal B}\coth(t{\cal B})\,
\label{337xx}
\ee
and $\left<u,\gamma v\right>=\gamma_{AB}u^Av^B$ is the inner product on
the algebra ${\cal G}$.

\begin{theorem}
Let $\Phi(t;k)$ be a function on the group $G$ defined by
\bea
\Phi(t;k)&=&(4\pi t)^{-N/2}A(t;k)
\exp\left(-\frac{\Theta(t;k)}{2t}+\frac{1}{6}R_G t\right)\,,
\label{46a}
\eea
Then $\Phi(t;k)$ satisfies the equation
\be
\partial_t \Phi = J^2\Phi\,,
\label{423}
\ee
and the initial condition
\be
\Phi(0;k)=|\gamma|^{-1/2}\delta(k)\,.
\label{428}
\ee
\end{theorem}

\noindent
\noindent\noindent\noindent{\bf Proof.}
We compute first
\be
\partial_t\Theta=\frac{1}{t}\Theta
-\frac{1}{2t}\left<k,\gamma\hat\Theta^2 k\right>
+\frac{t}{2}\left<k,\gamma{\cal B}^2 k\right>
\,
\ee
and
\be
\partial_t A=
\frac{1}{2t}\left(N
-\tr_{\cal G}\hat\Theta\right)A
\,.
\ee
Therefore,
\be
\partial_t\Phi=
\left[
\frac{1}{6}R_G-\frac{1}{2t}\tr_{\cal G}\hat\Theta
+\frac{1}{4t^2}\left<k,\gamma\hat\Theta^2k\right>
-\frac{1}{4}\left<k,\gamma{\cal B}^2k\right>
\right]\Phi\,.
\label{339xx}
\ee

Next, we have
\be
J^2=X^2-\gamma^{AB}{\cal B}_{AC}k^C X_B
+\frac{1}{4}\gamma^{AB}{\cal B}_{AC}{\cal B}_{BD}k^Ck^D\,.
\ee
By using the eqs. (\ref{327xx}) and  (\ref{2172})
and the anti-symmetry of the
matrix ${\cal B}_{AB}$ we show that
\be
\gamma^{AB}
{\cal B}_{AC}k^C X_B\Theta=0\,,
\ee
and
\be
\gamma^{AB}
{\cal B}_{AC}k^C X_B A=0\,,
\ee
and, therefore,
\be
{\cal B}_{AC}k^C X_B\Phi=0\,.
\ee
Thus,
\bea
J^2\Phi&=&
\Biggl[
A^{-1}(X^2A)
-\frac{1}{2t}(X^2\Theta)
+\frac{1}{4t^2}\gamma^{AB}(X_A\Theta)(X_B\Theta)
\nonumber\\
&&
-\frac{1}{t}A^{-1}\gamma^{AB}(X_B A) (X_A \Theta)
-\frac{1}{4}\left<k,\gamma{\cal B}^2k\right>
\Biggr]\Phi\,.
\label{344xx}
\eea

Further, by using \ref{327xx}) we get
\be
\gamma^{AB}(X_A\Theta)(X_B\Theta)
=\left<k,\gamma\hat\Theta^2 k\right>\,
\label{345xx}
\ee
\be
X^2\Theta=\tr_{\cal G}X+\tr_{\cal G}\hat\Theta-N\,.
\ee
Now, by using the eq. (\ref{313xx}) in the form
\be
A^2\partial_M(A^{-2}X_B{}^M)=0
\ee
and eqs. (\ref{2172}) and (\ref{315xx}) we show that
\be
A^{-1}\gamma^{AB}(X_A\Theta)X_B A=\frac{1}{2}\left(
N-\tr_{\cal G}X\right)\,.
\ee
Finally, by using eq. (\ref{425}) we obtain
\be
A^{-1}X^2 A =\frac{1}{6}R_G\,.
\label{349xx}
\ee

Finally, substituting the eqs. (\ref{345xx})-(\ref{349xx}) into eq.
(\ref{344xx}) and comparing it with eq. (\ref{339xx}) we prove the
eq. (\ref{423}). The initial condition (\ref{428}) follows easily from
the well known property of the Gaussian.
This completes the proof of the theorem.

%===========================================================
\subsection{Regularization and Analytical Continuation}

In the following we will complexify the gauged  curvature group in the
following sense. We extend the  canonical coordinates
$(k^A)=(p^a,\omega^i)$ to the whole complex Euclidean space  $\CC^N$.
Then all group-theoretic functions introduced above become analytic
functions of $k^A$ possibly with some poles on the real section $\RR^N$
for compact groups.  In fact, we replace the actual real slice $\RR^N$
of $\CC^N$ with an $N$-dimensional subspace $\RR^N_{\rm reg}$ in $\CC^N$
obtained by rotating the real section $\RR^N$ counterclockwise in
$\CC^N$ by $\pi/4$. That is, we replace each coordinate $k^A$ by
$e^{i\pi/4}k^A$. In the complex domain the group becomes non-compact. We
call this procedure the {\it decompactification}. If the group is
compact, or has a compact subgroup, then this plane will cover the
original group infinitely many times. 

Since the metric $(\gamma_{AB})=\diag(\delta_{ab},\beta_{ij})$ is not
necessarily positive definite,  (actually, only the metric of the
holonomy group $\beta_{ij}$ is non-definite) we analytically continue
the function $\Phi(t;k)$ in the  complex plane of $t$ with a cut along
the negative imaginary axis so that
$-\pi/2<\arg\,t<3\pi/2$.  Thus, the function $\Phi(t;k)$ defines an
analytic function of $t$ and $k^A$. For the purpose of the following
exposition we shall consider $t$ to be {\it real negative}, $t<0$. This
is needed in order to make all integrals convergent and well defined and
to be able to do the analytical continuation.

As we will show below, the singularities occur only in the holonomy
group. This means that there is no need to complexify the coordinates
$p^a$. Thus, in the following we assume the coordinates $p^a$ to be real
and the coordinates $\omega^i$ to be complex, more precisely, to take
values in the $p$-dimensional subspace $\RR^p_{\rm reg}$ of $\CC^p$ 
obtained by rotating $\RR^p$ counterclockwise by $\pi/4$ in $\CC^p$ That
is, we have $\RR^N_{\rm reg}=\RR^n\times \RR^p_{\rm reg}$.

This procedure (that we call a regularization) with the nonstandard
contour of integration is necessary  for the convergence of the
integrals below since we are treating both the compact and the
non-compact symmetric spaces simultaneously. Remember, that, in general,
the nondegenerate diagonal matrix $\beta_{ij}$ is not positive definite.
The space $\RR^p_{\rm reg}$ is chosen in such a way to make the Gaussian
exponent purely imaginary. Then the indefiniteness of the matrix $\beta$
does not cause any problems. Moreover, the integrand does not have any
singularities on these contours. The convergence of the integral is
guaranteed by the exponential growth of the sine for imaginary argument.
These integrals can be computed then in the following way. The
coordinates $\omega^j$ corresponding to the compact directions  are
rotated further by another $\pi/4$ to imaginary axis and the coordinates
$\omega^j$ corresponding to the non-compact directions are rotated back
to the real axis. Then, for $t<0$ all the integrals below are well
defined and convergent and define an analytic function of $t$ in a complex
plane with a cut along the negative imaginary axis.

%============================================================
\subsection{Heat Semigroup}

\begin{theorem}
The heat semigroup $\exp(t{\cal L}^2)$
can be represented 
in form of the integral
\be
\exp(t{\cal L}^2) = 
\int\limits_{\RR^N_{\rm reg}} dk\; |G|^{1/2}(k)\Phi(t;k)
\exp[{\cal L}(k)]\,.
\label{49a}
\ee
\end{theorem}

\noindent\noindent\noindent{\bf Proof.} Let
\be
\Psi(t) = 
\int\limits_{\RR^N_{\rm reg}} dk\; |G|^{1/2}\Phi(t;k)
\exp[{\cal L}(k)]\,.
\label{490}
\ee
By using the previous theorem we obtain
\be
\partial_t\Psi(t)
=\int\limits_{\RR^N_{\rm reg}} dk\; |G|^{1/2}\exp[{\cal L}(k)]
J^2\Phi(t;k)\,.
\label{421}
\ee
Now, by integrating by parts we get
\be
\partial_t\Psi(t)
=\int\limits_{\RR^N_{\rm reg}} dk\; |G|^{1/2}\Phi(t;k)
J^2\exp[{\cal L}(k)]\,,
\label{421a}
\ee
and, by using eq. (\ref{330xx}) we obtain
\be
\partial_t\Psi(t)=\Psi(t){\cal L}^2\,.
\ee
Finally from the initial condition (\ref{428}) for the function
$\Phi(t;k)$ we get
\be
\Psi(0)=1\,,
\ee
and, therefore, $\Psi(t)=\exp(t{\cal L}^2)$.

\begin{theorem}
Let $\Delta$ be the Laplacian acting on sections of a homogeneous
twisted spin-tensor vector bundle over a symmetric space. Then the
heat semigroup $\exp(t\Delta)$
can be represented in form of an integral
\bea
\exp(t\Delta) 
&=& 
(4\pi t)^{-N/2}
\det{}_{TM}\left(\frac{\sinh(t{\cal B})}{t{\cal B}}\right)^{-1/2} 
\exp\left(-t{\cal R}^2+ {1\over 6} R_G t\right)
\nonumber\\
&&
\int\limits_{\RR^N_{\rm reg}} dk\; |\gamma|^{1/2}
\det{}_{\cal G}\left({\sinh[C(k)/2]\over C(k)/2}\right)^{1/2}
\nonumber\\
&&
\times
\exp\left\{ -{1\over 4t}
\left<k,\gamma t{\cal B}\coth(t{\cal B})k\right>
\right\}
\exp[{\cal L}(k)]\,.
\label{49}
\eea
\end{theorem}

\noindent
\noindent\noindent{\bf Proof.}
By using the eq. (\ref{2110}) we obtain
\be
\exp(t\Delta)=\exp\left(-t{\cal R}^2\right)
\exp\left(t{\cal L}^2\right)\,.
\ee
The statement of the theorem follows now from the eqs. (\ref{49a}),
(\ref{46a}),  (\ref{334xx})-(\ref{337xx}) and (\ref{419}).

%=============================================================
\section{Heat Kernel}
\setcounter{equation}0

\subsection{Heat Kernel Diagonal and Heat Trace}

The heat kernel diagonal on a homogeneous bundle over a symmetric space
is parallel. In a local parallel local frame it is just a constant 
matrix. The fiber trace of the heat kernel diagonal is just a constant.
That is why, it can be computed at any point in $M$. We fix a point $x'$
in $M$ such that the Killing vectors satisfy the initial conditions 
(\ref{237})-(\ref{238}) and are given by the explicit formulas above
(\ref{226})-(\ref{229}). We compute the heat kernel diagonal at the
point $x'$.

The heat kernel diagonal can be obtained by acting by the heat
semigroup $\exp(t\Delta)$ on the delta-function,
\cite{avramidi94,avramidi96}
\bea
U^{\rm diag}(t)&=&\exp(t\Delta)\delta(x,x')\Big|_{x=x'}
\nonumber\\
&=&
\exp\left(-t{\cal R}^2\right)
\int\limits_{\RR^N_{\rm reg}} dk\; |G|^{1/2}\Phi(t;k)
\exp[{\cal L}(k)]\delta(x,x')\Big|_{x=x'}\,.
\label{41}
\eea
To be able to use this integral representation we need to compute the
action of the isometries $\exp[{\cal L}(k)]$ on the delta-function.

\begin{proposition}
Let $\varphi$ be a section of the twisted spin-tensor bundle 
${\cal V}$, ${\cal L}_A$ be the twisted Lie derivatives,
$k^A=(p^a,\omega^i)$ be the canonical coordinates on the group
and
${\cal L}(k)=k^A{\cal L}_A$.
Let $\xi=k^A\xi_A$ be the Killing vector and $\psi_t$
be the corresponding one-parameter diffeomorphism.
Then
\bea
\exp\left[{\cal L}(k)\right]\varphi(x)
&=&\exp\left(-\frac{1}{2}\theta_{ab}G^{ab}\right)
\varphi(\hat x)\Bigg|_{t=1}\,,
\label{2135}
\eea
where $\hat x=\psi_t(x)$ and the matrix $\theta$ is defined by
(\ref{2234mm}).
In particular, for $p=0$ and $x=x'$
\bea
\exp[{\cal L}(k)]\varphi(x)\Big|_{p=0,x=x'}
&=&
\exp\left[{\cal R}(\omega)\right]
\varphi(x)\,.
\label{2150}
\eea
\end{proposition}

\noindent\noindent\noindent{\bf Proof.}
This statement follows from eqs. (\ref{2236mm}) and (\ref{2237mm})
and the fact that the Lie derivative is nothing 
but the generator of the pullback.

\begin{proposition}
Let $\omega^i$ be the canonical coordinates on the holonomy
group $H$ and $(k^A)=(p^a,\omega^i)$ be the natural splitting of the
canonical coordinates on the curvature group $G$.
Then
\be
\exp[{\cal L}(k)]\delta(x,x')\Big|_{x=x'}
=\det{}_{TM}\left(\frac{\sinh[D(\omega)/2]}{D(\omega)/2}\right)^{-1}
\exp[{\cal R}(\omega)]
\delta(p)\,.
\label{424a}
\ee
\end{proposition}

\noindent
\noindent\noindent\noindent{\bf Proof.}
Let $\hat x(t,p,\omega,x,x')=\psi_t(x)$.
By making use of the eq. (\ref{2135})
we obtain
\bea
\exp[{\cal L}(k)]\delta(x,x')\Big|_{x=x'}
&=&\exp\left(-\frac{1}{2}\theta_{ab}G^{ab}\right)
\delta(\hat x(1,p,\omega,x,x'),x')\Big|_{x=x',t=1}\,.
\eea
Now we change the variables from $x^\mu$ 
to the normal coordinates $y^a$
to get
\be
\delta(\hat x(1,p,\omega,x,x'),x')\Big|_{x=x'}
=|g|^{-1/2}
\det\left(\frac{\partial y^a}{\partial x^\mu}\right)
\delta(\hat y(1,p,\omega,y))\Big|_{y=0}\,.
\ee
This delta-function picks the values of $p$ that make
$\hat y=0$, which is exactly the functions
$\bar p=\bar p(\omega,y)$
defined by the eq. (\ref{2197nn}). By switching further
to the  variables $p$ we obtain
\be
\delta(\hat x(1,p,\omega,x,x'),x')\Big|_{x=x'}
=|g|^{-1/2}
\det\left(\frac{\partial y^a}{\partial x^\mu}\right)
\det\left(\frac{\partial \hat y^b}{\partial p^c}\right)^{-1}
\delta(p-\bar p(\omega,y))\Big|_{y=0,t=1}\,.
\ee
Now, by recalling from (\ref{2216mm}) that 
$\bar p|_{y=0}=0$ and by using (\ref{244mm})
and (\ref{429a})
we evaluate the Jacobians for $p=y=0$ and $t=1$ to 
get the eq. (\ref{424a}).

{\bf Remarks.}
Some remarks are in order here. 
We implicitly assumed that there are no closed geodesics
and that the equation of closed orbits of isometries
\be
\hat y^a(1,\bar p,\omega,0)=0\,
\label{410}
\ee
has a unique solution $\bar p=\bar p(\omega,0)=0$. On compact symmetric
spaces this is not true: there are infinitely many closed geodesics and
infinitely many closed orbits of isometries.   However, these global
solutions, which reflect the global topological structure of the
manifold, will not affect our local analysis.  In particular, they do
not affect the asymptotics of the heat kernel.  That is why, we have
neglected them here.  This is reflected in the fact that the Jacobian in
(\ref{424a}) can become singular when the coordinates of the holonomy
group $\omega^i$ vary from $-\infty$ to $\infty$.  Note that the exact
results for compact symmetric spaces can be obtained by an analytic
continuation from the dual noncompact case when such closed geodesics
are absent \cite{camporesi90}. That is why we proposed above to
complexify our holonomy group. If the coordinates $\omega^i$ are complex
taking values in the subspace $\RR^p_{\rm reg}$ defined above, then the
equation (\ref{410}) should have a unique solution and the Jacobian is
an analytic function. It is worth stressing once again  that the
canonical coordinates cover the whole group except for a set of measure
zero.  Also a compact subgroup is covered  infinitely many times. 
We will show below how this works in the case of the two-sphere,
$S^2$.

Now by using the above lemmas and the theorem we can compute the heat
kernel diagonal. 
We define the
matrix  $F(\omega)$ by
\be
F(\omega)=\omega^i F_i\,.
\ee

\begin{theorem}
The heat kernel diagonal of the 
Laplacian on twisted spin-vector bundles over a symmetric space has the form
\bea
U^{\rm diag}(t)
&=&(4\pi t)^{-n/2}
\det{}_{TM}\left(\frac{\sinh(t{\cal B})}{t{\cal B}}\right)^{-1/2}
\exp\left\{\left({1\over 8} R + {1\over 6} R_H 
-{\cal R}^2\right)t\right\}
\nonumber
\label{437a}
\\
&&
\times
\int\limits_{\RR^n_{\rm reg}} 
\frac{d\omega}{(4\pi t)^{p/2}}\;|\beta|^{1/2}
\exp\left\{-{1\over 4 t}\left<\omega,\beta\omega\right>\right\}
\cosh\left[\,{\cal R}(\omega)\right]
\nonumber\\
&& 
\times\det{}_{\cal H}
\left({\sinh\left[\,F(\omega)/2\right]\over 
\,F(\omega)/2}\right)^{1/2}
\det{}_{TM}\left({\sinh\left[\,D(\omega)/2\right]\over 
\,D(\omega)/2}\right)^{-1/2}\,,
%\nonumber
\eea
where $|\beta|=\det \beta_{ij}$ and
$\left<\omega,\beta\omega\right>=\beta_{ij}\omega^i\omega^j$.
\end{theorem}

\noindent
\noindent\noindent\noindent{\bf Proof.}
First, we have $dk=dp\;d\omega$ and
$|\gamma|=|\beta|\,.$
By using the equations (\ref{41}) 
and (\ref{424a}) and integrating over $p$ we obtain 
the heat kernel diagonal
\be
U^{\rm diag}(t)
=\int\limits_{\RR^p_{\rm reg}}d\omega\;|G|^{1/2}(0,\omega)\Phi(t;0,\omega)
\det{}_{TM}\left(\frac{\sinh[D(\omega)/2]}{D(\omega)/2}\right)^{-1}
\exp[{\cal R}(\omega)-t{\cal R}^2]\,.
\ee
Further, by using the eq. (\ref{318}) 
we compute the determinants
\be
\det{}_{\cal G}\left({\sinh[C(\omega)/2]\over C(\omega)/2}\right)
=\det{}_{TM}\left({\sinh[D(\omega)/2]\over D(\omega)/2}\right)
\det{}_{\cal H}\left({\sinh[F(\omega)/2]\over F(\omega)/2}\right)\,.
\label{437aa}
\ee
Now, we by using (\ref{2171xx}) we compute (\ref{337mm}) 
\be
\Theta(t;0,\omega)=\frac{1}{2}\left<\omega,\beta\omega\right>\,,
\ee
and, finally, 
by using eq. (\ref{46a}), (\ref{334xx}), (\ref{312xx})
and (\ref{295vv}) 
we get the result (\ref{437a}).

By using this theorem we can also compute the heat trace
for compact manifolds
\bea
\Tr_{L^2}\exp(t\Delta)
&=& \int\limits_M d\vol\;
(4\pi t)^{-n/2}\tr_V
\det{}_{TM}\left(\frac{\sinh(t{\cal B})}{t{\cal B}}\right)^{-1/2}
%\nonumber
\\
&&
\times
\exp\left\{\left({1\over 8} R + {1\over 6} R_H 
-{\cal R}^2
\right)t\right\}
\nonumber\\
&&
\times
\int\limits_{\RR^p_{\rm reg}} 
\frac{d\omega}{(4\pi t)^{p/2}}\;|\beta|^{1/2}
\exp\left\{-{1\over 4 t}\left<\omega,\beta\omega\right>\right\}
\cosh\left[{\cal R}(\omega)\right]
\nonumber\\
&& 
\times\det{}_{\cal H}
\left({\sinh\left[F(\omega)/2\right]\over 
F(\omega)/2}\right)^{1/2}
\det{}_{TM}\left({\sinh\left[D(\omega)/2\right]\over 
D(\omega)/2}\right)^{-1/2}
\,,
\label{437cc}
\nonumber
\eea
where $\tr_V$ is the fiber trace.

%=====================================================
\subsection{Heat Kernel Asymptotics}

It is well known that there is the following asymptotic
expansion as $t\to 0$ of the heat kernel diagonal
\cite{gilkey95}
\be
U^{\rm diag}(t)\sim (4\pi t)^{-n/2}
\sum_{k=0}^\infty t^k a_k\,.
\ee
The coefficients $a_k$ are called the local heat kernel
coefficients.
On compact manifolds, there is a similar
asymptotic expansion of the heat trace
with the global heat invariants $A_k$ defined by
\be
A_k=\int_M d\vol\; \tr_V a_k\,.
\ee
In symmetric spaces the heat invariants do not contain any
additional information since the local heat kernel coefficients
define the heat invariants $A_k$ up to a constant equal
to the volume of the manifold, 
\be
A_k=\vol(M)\tr_V a_k\,.
\ee

We introduce a Gaussian average over the holonomy algebra by
\be
\left<f(\omega)\right> = \int\limits_{\RR^p_{\rm reg}}
\frac{d\omega}{(4\pi)^{p/2}}\; |\beta|^{1/2}
\exp\left(-{1\over 4}\left<\omega,\beta\omega\right>\right)f(\omega)
\ee
Then we can write
\bea
&&U^{\rm diag}(t) 
=(4\pi t)^{-n/2}  
\det{}_{TM}\left(\frac{\sinh(t{\cal B})}{t{\cal B}}\right)^{-1/2}
\exp\left\{\left({1\over 8} R + {1\over 6} R_H 
-{\cal R}^2\right) t \right\}
\label{422}
\\[12pt]
&&
%\hspace{-1cm}
\times
\Bigg< \cosh\left[\sqrt{t}\,{\cal R}(\omega)\right]
\det{}_{\cal H}\left({\sinh\left[\sqrt{t}\,F(\omega)/2\right]\over 
\sqrt{t}\,F(\omega)/2}\right)^{1/2}
\det{}_{TM}\left({\sinh\left[\sqrt{t}\,D(\omega)/2\right]\over 
\sqrt{t}\,D(\omega)/2}\right)^{-1/2}\Bigg>   
\nonumber
\eea
This equation can be used now to generate all heat kernel
coefficients $a_k$ for any locally symmetric space  simply by expanding
it in a power  series in $t$. By using the standard Gaussian averages
\footnote{We have corrected here a misprint in the eq. (4.68) of 
\cite{avramidi96}.}
\bea
\left<\omega^i_1\cdots \omega^{i_{2k+1}}\right> &=& 0\,,
\label{440a}
\\
\left<\omega^{i_1}\cdots \omega^{i_{2k}}\right> 
&=& {(2k)!\over k!}\beta^{(i_1 i_2}\cdots
\beta^{i_{2k-1}i_{2k})}\,,
\label{441a}
\eea
one can obtain now all heat kernel coefficients in terms of 
traces of various contractions 
of the matrices $D^a{}_{ib}$ and $F^j{}_{ik}$ with  the
matrix  $\beta^{ik}$.
All these quantities are curvature invariants and
can be  expressed directly in terms of the Riemann tensor.

There is an alternative representation of the Gaussian average
in purely algebraic terms.
Let $b^j$ and $b^*_k$ be operators, 
called creation and annihilation operators,
acting on a Hilbert space,
that satisfy the following commutation relations
\be
[b^j,b^*_k]=\delta^j_k\,,
\ee
\be
[b^j,b^k]=[b^*_j,b^*_k]=0\,.
\ee
Let $|0\rangle$ be a unit vector in the Hilbert space, called the vacuum
vector, that satisfies the equations
\be
\langle 0|0\rangle=1\,,
\ee
\be
b^j|0\rangle=\langle 0|b_k^*=0\,.
\ee
Then the Gaussian average is nothing but the vacuum expectation
value 
\be
\langle f(\omega)\rangle
=\langle 0| f(b)\exp\langle b^*,\beta b^*\rangle |0\rangle\,,
\ee
where $\langle b^*,\beta b^*\rangle=\beta^{jk}b^*_jb^*_k$\,.
This should be computed by the so-called normal ordering, that is, 
by simply commuting the operators $b_j$
through the operators $b^*_k$ until they hit the vacuum vector
giving zero. The remaining non-zero commutation terms precisely
reproduce the eqs. (\ref{440a}), (\ref{441a}).

%===============================================
\subsubsection{Calculation of the Coefficient $a_1$}

As an example let us calculate the lowest heat kernel coefficients:
$a_0$ and $a_1$.
Let $X$ be a matrix. Then by using
\be
\det\left({\sinh \left(\sqrt{t}\,X\right)\over \sqrt{t}\,X}\right)^{m}
=\exp\left(m\;\tr\log {\sinh \left(\sqrt{t}\,X\right)\over \sqrt{t}\,X}\right)
\ee
and \cite{erdelyi53}
\be
\log{\sinh\left(\sqrt{t}\,X\right)\over \sqrt{t}\,X}
=\sum_{k=1}^\infty \frac{2^{2k-1}B_{2k}}{k(2k)!} t^k X^{2k}\,,
\ee
where $B_k$ are Bernoulli numbers, in particular,
\be
B_0=1,\qquad
B_1=-\frac{1}{2},\qquad
B_2=\frac{1}{6},
\ee
we obtain
\bea
\det\left({\sinh \left(\sqrt{t}\,X\right)\over 
\sqrt{t}\,X}\right)^{\pm 1/2}
&=&
1\pm\frac{1}{12}t\, \tr X^{2}
+O(t^2)\,.
\eea
Now, by using eq. (\ref{422}) we obtain
\be
U^{\rm diag}(t)
=(4\pi t)^{-n/2}\left[a_0+t a_1+O(t^2)\right]\,,
\ee
where
\be
a_0=\II\,,\qquad
a_1=\left<b_1\right>\,,
\ee
and
\bea
b_1=
\left[{1\over 8} R +{1\over 6} R_H
+\frac{1}{48}\tr F(\omega)^{2}
-\frac{1}{48}\tr D(\omega)^{2}
\right]\II
-{\cal R}^2
+\frac{1}{2}\;{\cal R}(\omega)^2
\,.
\label{437bb}
\eea
Next, bu using (\ref{441a}), in particular,
\be
\left<\omega^i\omega^j\right>=2\beta^{ij},
\ee
we obtain
\be
\left<{\cal R}(\omega)^2\right>=2{\cal R}^2\,,
\ee
\be
\left<\tr F(\omega)^2\right>=2\tr F_i F^i
=2F^j{}_{il}F^{li}{}_j=-8R_H\,,
\ee
\be
\left<\tr D(\omega)^2\right>=2\tr D_i D^i
=2 D^a{}_{ib}D^{bi}{}_a=-2R\,,
\ee
and, therefore,
\bea
a_1&=&
\left[{1\over 8} R+{1\over 6} R_H
-\frac{1}{6}R_H
+\frac{1}{24}R\right]\II
-{\cal R}^2+{\cal R}^2
%\nonumber\\
%&=&
=\frac{1}{6}R\II\,.
\eea
This confirms the well know result for the coefficient $a_1$
\cite{gilkey95,avramidi91}.

%=============================================================
\subsection{Heat Kernel on $S^2$ and $H^2$}

Let us apply our result to a special case of a two-sphere $S^2$ of
radius $r$, which is  a compact
symmetric space equal to the quotient of the
isometry group, $SO(3)$, by the isotropy group, $SO(2)$,
\be
S^2=SO(3)/SO(2).
\ee
The two-sphere is too small to incorporate an additional
Abelian field ${\cal B}$; therefore, we set ${\cal B}=0$.

Let $y^a$ be the normal coordinates defined above.
On the 2-sphere of radius $r$ they range over 
$-r\pi\le y^a\le r\pi$. We define the polar coordinates 
$\rho$ and $\varphi$
by
\be
y^1=\rho\cos\varphi,\qquad
y^2=\rho\sin\varphi\,,
\ee
so that $0\le \rho\le r\pi $ and $0\le \varphi\le 2\pi$.

The orthonormal frame of $1$-forms is
\be
e^1=d\rho\,,\qquad
e^2=r\sin\left(\frac{\rho}{r}\right)d\varphi\,,
\ee
which gives the spin connection $1$-form
\be
\omega_{ab}=-\varepsilon_{ab}
\cos\left(\frac{\rho}{r}\right)d\varphi\,
\,
\ee
with $\varepsilon_{ab}$ being the antisymmetric Levi-Civita tensor,
and the curvature
\be
R_{abcd}=\frac{1}{r^2}\varepsilon_{ab}\varepsilon_{cd}
=\frac{1}{r^2}(\delta_{ac}\delta_{bd}-\delta_{ad}\delta_{bc})\,,
\ee
\be
R_{ab}=\frac{1}{r^2}\delta_{ab}\,,
\qquad
R=\frac{2}{r^2}\,.
\ee

Since the holonomy group $SO(2)$ is one-dimensional, it is obviously
Abelian, so all structure constants $F^i{}_{jk}$ are equal to zero, and
therefore, the curvature of the holonomy group vanishes, $R_H=0$. The
metric of the holonomy group $\beta_{ij}$ is now just a constant,
$\beta=1/r^2$.  The only generator of the holonomy group in the
vector representation is
\be
D_{ab}=-\frac{1}{r^2}E_{ab}=-\frac{1}{r^2}\varepsilon_{ab}\,.
\ee

The irreducible representations of $SO(2)$ are parametrized by
$\alpha$, which is either an integer, $\alpha=m$, or a half-integer,
$\alpha=m+\frac{1}{2}$. 
Therefore, the generator ${\cal R}$ of the holonomy group
and the Casimir operator ${\cal R}^2$ are
\be
{\cal R}=i\frac{\alpha}{r^2}\,,
\ee
\be
{\cal R}^2=\beta^{ij}{\cal R}_i{\cal R}_j=-\frac{\alpha^2}{r^2}\,.
\ee
The extra factor $r^2$ here is due to the inverse metric
$\beta^{-1}=r^2$ of the holonomy group.

The Lie derivatives ${\cal L}_A$  are given by 
\be 
{\cal L}_1=\cos\varphi\partial_\rho
-\frac{\sin\varphi}{r}\cot\left(\frac{\rho}{r}\right)
\partial_\varphi
+i\frac{\sin\varphi}{r\sin\left(\rho/r\right)}\alpha
\,, 
\ee 
\be 
{\cal L}_2=\sin\varphi\partial_\rho
+\frac{\cos\varphi}{r}\cot\left(\frac{\rho}{r}\right)
\partial_\varphi
-i\frac{\cos\varphi}{r\sin\left(\rho/r\right)}\alpha
\,,
\ee 
\be 
{\cal L}_3=\frac{1}{r^2}\partial_\varphi
\,, 
\ee 
and form a
representation of the $SO(3)$ algebra 
\be 
[{\cal L}_1,{\cal L}_2]=-{\cal
L}_3\,, \qquad [{\cal L}_3,{\cal L}_1]=-\frac{1}{r^2}{\cal L}_2 
\qquad
[{\cal L}_3,{\cal L}_2]=\frac{1}{r^2}{\cal L}_1 \,. 
\ee
The Laplacian is
given by 
\be 
\Delta=\partial_\rho^2
+\frac{1}{r}\cot\left(\frac{\rho}{r}\right)\partial_\rho 
+\frac{1}{r^2\sin^2\left(\rho/r\right)}
\left[\partial_\varphi
-i\alpha\cos\left(\frac{\rho}{r}\right)\right]^2
\ee 

Now, we need to compute the determinant
\be
\det{}_{TM}\left(\frac{\sinh[\omega D]}{\omega D}\right)^{-1/2}
=\frac{\omega/(2r^2)]}{\sin[\omega/(2r^2)]}\,.
\ee
The contour of integration  over $\omega$ in (\ref{437a}) should be the
real axis rotated counterclockwise by $\pi/4$. Since $S^2$ is compact,
we rotate it further to the imaginary axis and rescale $\omega$ for
$t<0$ by $\omega\to r\sqrt{-t}\,\omega$ to obtain an analytic function
of $t$
\bea
U^{\rm diag}(t)&=&\frac{1}{4\pi t}
\exp\left[\left(\frac{1}{4}+\alpha^2\right)\frac{t}{r^2}\right]
%\nonumber
\\
&&\times
\int\limits\limits_{-\infty}^\infty
\frac{d\omega}{\sqrt{4\pi}}
\exp\left(-\frac{\omega^2}{4}\right)
\frac{\omega\sqrt{-t}/(2r)}
{\sinh\left[\omega\sqrt{-t}/(2r)\right]}
\cosh\left(\alpha\omega\sqrt{-t}/r\right)
\,.
\nonumber
\eea
If we would have rotated the contour to the real axis instead then we
would have obtained after rescaling 
$\omega\to r\sqrt{t}\,\omega$
for $t>0$,
\bea
U^{\rm diag}(t)&=&\frac{1}{4\pi t}
\exp\left[\left(\frac{1}{4}+\alpha^2\right)\frac{t}{r^2}\right]
%\nonumber
\\
&&\times
\fint\limits\limits_{-\infty}^\infty
\frac{d\omega}{\sqrt{4\pi}}
\exp\left(-\frac{\omega^2}{4}\right)
\frac{\omega \sqrt{t}/(2r)}
{\sin\left[\omega \sqrt{t}/(2r)\right]}
\cos\left(\alpha\omega\sqrt{t}/r\right)
\,,
\nonumber
\eea
where $\fint$ denotes the Cauchy principal value of the integral. This
can also be written as
\bea
&&
U^{\rm diag}(t)=\frac{1}{4\pi t}
\exp\left[\left(\frac{1}{4}+\alpha^2\right)\frac{t}{r^2}\right]
%\nonumber
\\
&&
\times
\sum_{k=-\infty}^\infty (-1)^k
\int\limits\limits_{0}^{2\pi r/\sqrt{t}}
\frac{d\omega}{\sqrt{4\pi}}
\exp\left[-\frac{1}{4}
\left(\omega+\frac{2\pi r}{\sqrt{t}}k\right)^2\right]
\frac{\sqrt{t}}{2r}
\frac{\left(\omega+\frac{2\pi r}{\sqrt{t}}k\right)}
{\sin\left[\omega\sqrt{t}/(2r)\right]}
\nonumber\\
&&
\times
\cos\left(\alpha\omega\sqrt{t}/r\right)
\,.
\nonumber
\eea
This is nothing but the sum over the closed geodesics of $S^2$.
Note that the factor $\cos\left(\alpha\omega\sqrt{t}/r\right)$
is either periodic (for integer $\alpha$) or anti-periodic
(for half-integer $\alpha$).

%===========================
The non-compact symmetric space dual to the $2$-sphere is the
hyperbolic plane $H^2$ of pseudo-radius $a$. 
It is equal to the quotient of the
isometry group, $SO(1,2)$, by the isotropy group, $SO(2)$,
\be
H^2=SO(1,2)/SO(2).
\ee
Let $y^a$ be the normal coordinates defined above.
On $H^2$ they range over 
$-\infty\le y^a\le \infty$. We define the polar coordinates 
$u$ and $\varphi$
by
\be
y^1=u\cos\varphi,\qquad
y^2=u\sin\varphi\,,
\ee
so that $0\le u\le \infty $ and $0\le \varphi\le 2\pi$.

The orthonormal frame of $1$-forms is
\be
e^1=du\,,\qquad
e^2=a\sinh\left(\frac{u}{a}\right)d\varphi\,,
\ee
which gives the spin connection $1$-form
\be
\omega_{ab}=-\varepsilon_{ab}
\cosh\left(\frac{u}{a}\right)d\varphi\,
\,,
\ee
and the curvature
\be
R_{abcd}=-\frac{1}{a^2}\varepsilon_{ab}\varepsilon_{cd}
=-\frac{1}{a^2}(\delta_{ac}\delta_{bd}-\delta_{ad}\delta_{bc})\,,
\ee
\be
R_{ab}=-\frac{1}{a^2}\delta_{ab}\,,
\qquad
R=-\frac{2}{a^2}\,.
\ee
The metric of the isotropy group $\beta_{ij}$ is just a constant,
$\beta=-1/a^2$, and the only generator of the
isotropy group in the vector representation is given by
\be
D_{ab}=\frac{1}{a^2}E_{ab}
=\frac{1}{a^2}\varepsilon_{ab}\,.
\ee

The Lie derivatives ${\cal L}_A$ are now
\be
{\cal L}_1=\cos\varphi\partial_u
-\frac{\sin\varphi}{a}\coth\left(\frac{u}{a}\right)
\partial_\varphi
+i\frac{\sin\varphi}{a\sinh\left(u/a\right)}\alpha
\,,
\ee
\be
{\cal L}_2=\sin\varphi\partial_u
+\frac{\cos\varphi}{a}\coth\left(\frac{u}{a}\right)
\partial_\varphi
-i\frac{\cos\varphi}{a\sinh\left(u/a\right)}\alpha
\,,
\ee
\be
{\cal L}_3=-\frac{1}{a^2}\partial_\varphi
\,,
\ee
and form a representation of the $SO(1,2)$
algebra
\be
[{\cal L}_1,{\cal L}_2]=-{\cal L}_3\,,
\qquad
[{\cal L}_3,{\cal L}_1]=\frac{1}{a^2}{\cal L}_2
\qquad
[{\cal L}_3,{\cal L}_2]=-\frac{1}{a^2}{\cal L}_1
\,.
\ee
The Laplacian is given by
\be
\Delta=\partial_u^2
+\frac{1}{a}\coth\left(\frac{u}{a}\right)\partial_u
+\frac{1}{a^2\sinh^2\left(u/a\right)}
\left[\partial_\varphi
-i\alpha\cosh\left(\frac{u}{a}\right)\right]^2
\ee

The contour of integration over $\omega$ in (\ref{437a})
for the heat kernel should be the real axis
rotated counterclockwise  by $\pi/4$. Since $H^2$ is non-compact, we
rotate it back to the real axis and rescale $\omega$ for $t>0$ by
$\omega\to a\sqrt{t}\,\omega$ to obtain the heat kernel diagonal for the
Laplacian on $H^2$
\bea
U^{\rm diag}(t)&=&\frac{1}{4\pi t}
\exp\left[-\left(\frac{1}{4}+\alpha^2\right)\frac{t}{a^2}\right]
%\nonumber
\label{471cc}
\\
&&\times
\int\limits\limits_{-\infty}^\infty
\frac{d\omega}{\sqrt{4\pi}}
\exp\left(-\frac{\omega^2}{4}\right)
\frac{\omega\sqrt{t}/(2a)}
{\sinh\left[\omega\sqrt{t}/(2a)\right]}
\cosh\left(\alpha\omega\sqrt{t}/a\right)
\,.
\nonumber
\eea

We see that the heat kernel in the compact case of the two-sphere,
$S^2$, is related  with the heat  kernel in the non-compact case of the
hyperboloid, $H^2$, by the analytical continuation, $a^2\to -r^2$, or
$a\to ir$, or, alternatively, by replacing $t\to -t$ (and $a=r$).  One
can go even further and compute the Plancherel (or Harish-Chandra)
measure in the case of $H^2$ and the spectrum in the case of
$S^2$.

For $H^2$ we rescale the integration variable in (\ref{471cc}) by 
$\omega\to \omega a/\sqrt{t}$, substitute
\be
\frac{a}{\sqrt{4\pi t}}
\exp\left(-\frac{a^2 }{4t}\omega^2\right)
=\int\limits_{-\infty}^\infty
\frac{d\nu}{2\pi}\exp\left(-\frac{t}{a^2}\nu^2+i\omega\nu\right)\,,
\ee
integrate by parts over $\nu$,
and use
\bea
\int\limits\limits_{-\infty}^\infty
\frac{d\omega}{2\pi i}\;
e^{i\omega\nu}
\frac{\cosh(\alpha\omega)}
{\sinh\left(\omega/2\right)}
=\frac{1}{2}\left\{\tanh[\pi(\nu+i\alpha)]
+\tanh[\pi(\nu-i\alpha)]
\right\}\,
\eea
(and the fact that $\alpha$ is a half-integer)
to represent the heat kernel for $H^2$ in the form
\be
U^{\rm diag}(t)=\frac{1}{4\pi a^2}
\int\limits_{-\infty}^\infty
d\nu\;\mu(\nu)
\exp\left\{-\left(\frac{1}{4}+\alpha^2+\nu^2
\right)\frac{t}{a^2}
\right\}
\,,
\label{473cc}
\ee
where
\be
\mu(\nu)=\nu\tanh\nu
\ee
for integer $\alpha=m$, and
\be
\mu(\nu)=\nu\coth\nu
\ee
for half-integer $\alpha=m+\frac{1}{2}$.

For $S^2$ we proceed as follows. We cannot just substitute $a^2\to -r^2$
in (\ref{473cc}). Instead, first, we  deform the contour of integration
in (\ref{473cc}) to the $V$-shaped contour that consists of two segments
of straight lines, one going from $e^{i3\pi/4}\infty$ to $0$, and
another going from $0$ to $e^{i\pi/4}\infty$. Then, after we replace
$a^2\to -r^2$, we can  deform the contour further to go counterclockwise
around the positive imaginary axis. Then we notice that the function
$\mu(\nu)$ is a meromorphic function with simple poles on the imaginary
axis at $\nu_k=i d_k$, where
\be
d_k=\left(k+\frac{1}{2}\right), \qquad k=0,\pm 1\pm 2,\dots, \qquad
\mbox{ for integer } \alpha=m\,,
\ee
and at 
\be
d_k=k, \qquad k=\pm 1,\pm 2,\dots, \qquad  
\mbox{ for half-integer } \alpha=m+\frac{1}{2}\,.
\ee
Therefore, we can compute the integral by
residue theory to get
\be
U^{\rm diag}(t)=\frac{1}{4\pi r^2}
\sum_{k=0}^\infty d_k \exp\left(-\lambda_k t\right)\,,
\label{475cc}
\ee
where
\be
\lambda_k=\frac{1}{r^2}\left[\left(k+\frac{1}{2}\right)^2
-\frac{1}{4}-m^2\right]
\qquad
\mbox{ for integer } \alpha=m,
\ee
and
\be
\lambda_k=\frac{1}{r^2}\left[k^2
-\frac{1}{4}-\left(m+\frac{1}{2}\right)^2\right]
\qquad
\mbox{ for half-integer } \alpha=m+\frac{1}{2}\,.
\ee

Our results for the heat kernel on the 2-sphere $S^2$ and the hyperbolic
plane $H^2$ coincide with the exact  heat kernel of scalar Laplacian
(when ${\cal R}=\alpha=0$) reported in \cite{camporesi90} and obtained
by completely different methods.

%========================================================
%========================================================
\subsection{Index Theorem}

We can now apply this result for the calculation of the index of the
Dirac operator on spinors on compact manifolds
\be
D=\gamma^\mu\nabla_\mu\,.
\ee
Let the dimension $n$ of the manifold  be even and
\be
\Gamma=\frac{1}{n!} i^{n(n-1)/2}
\varepsilon^{a_1\dots a_n}
\gamma_{[a_1}\cdots\gamma_{a_n]}
\ee
be the chirality operator of the spinor representation so that
\be
\Gamma^2=\II_S
\ee
and
\be
\Gamma\gamma_a=-\gamma_a\Gamma\,.
\ee

Then the index of the Dirac operator is equal to
\be
\Ind(D)=\Tr_{L^2}\Gamma\exp(tD^2)\,.
\ee
We compute the square of the Dirac operator
by using the eqs. (\ref{213mm}), (\ref{219mm}), (\ref{117}) and (\ref{2153mm})
\bea
D^2&=&\Delta-\frac{1}{4}R\,\II_S
+\frac{1}{2}{\cal F}_{ab}\gamma^{ab}
\nonumber\\
&=&\Delta-\frac{1}{4}R\,\II_S
-\frac{1}{2}E^i{}_{ab}T_i\gamma^{ab}
+\frac{1}{2}\gamma^{ab}{\cal B}_{ab}\,.
\eea
In this case the generators ${\cal R}_i$ have the form
\be
{\cal R}_i=-\frac{1}{4}D^a{}_{ib}\gamma^b{}_a\otimes \II_T+\II_S\otimes T_i
\ee
and
the Casimir operator of the holonomy group in the spinor representation
is obtained by using (\ref{117})
\be
{\cal R}^2=\frac{1}{8}R\,\II_S
+\II_S\otimes T^2
-\frac{1}{2}E^j{}_{ab}\gamma^{ab}\otimes T_j
\,.
\ee

Thus, we obtain the index
\bea
\Ind(D)
&=& \int\limits_M d\vol\;
(4\pi t)^{-n/2}
\tr_V
\Gamma
\det{}_{TM}\left(\frac{\sinh(t{\cal B})}{t{\cal B}}\right)^{-1/2}
\nonumber\\
&&
\times
\exp\left\{\left(-{1\over 4} R + {1\over 6} R_H 
-T^2
+\frac{1}{2}{\cal B}_{ab}\gamma^{ab}
\right)t\right\}
\nonumber\\
&&
\times
\int\limits_{\RR^p_{\rm reg}} \frac{d\omega}{(4\pi t)^{p/2}}\;|\beta|^{1/2}
\exp\left\{-{1\over 4 t}\left<\omega,\beta\omega\right>\right\}
\nonumber\\
&& 
\times
\cosh\left(-\frac{1}{4}\omega^iD^a{}_{ib}\gamma^b{}_a
+\omega^i T_i\right)
\nonumber\\
&& 
\times\det{}_{\cal H}
\left({\sinh\left[F(\omega)/2\right]\over 
F(\omega)/2}\right)^{1/2}
\det{}_{TM}\left({\sinh\left[D(\omega)/2\right]\over 
D(\omega)/2}\right)^{-1/2}
\,.
\label{437d}
\eea

Since the index does not depend on $t$, the right-hand side of this
equation does not depend on $t$. By expanding it in an asymptrotic
power series in $t$, we see that the index is equal to
\be
\Ind(D)=(4\pi)^{-n/2}\int_M d\vol\, \tr_V\,\Gamma a_{n/2}\,.
\ee

%===================================================================
\section{Conclusion}

We have continued the study of the heat kernel on homogeneous spaces
initiated in 
\cite{avramidi93,avramidi94a,avramidi94,avramidi95,avramidi96}. In those
papers we have developed a systematic technique for calculation of the
heat kernel in two cases: a) a Laplacian on a vector bundle  with a
parallel curvature over a flat space \cite{avramidi93,avramidi95}, and
b) a scalar Laplacian on manifolds with parallel curvature
\cite{avramidi94,avramidi96}. What was missing in that study was the
case of  a non-scalar Laplacian on vector bundles with parallel
curvature over curved manifolds with parallel curvature. 

In the present paper we considered the Laplacian on a homogeneous bundle
and generalized the technique developed in \cite{avramidi96} to compute
the corresponding heat semigroup and the heat kernel. It is worth
pointing out that our formal result applies to general symmetric spaces
by making use of the regularization and the analytical continuation
procedure described above. Of course, the heat kernel coefficients are
just polynomials in the curvature and do not depend on this kind of
analytical continuation (for more detail, see \cite{avramidi96}).  

As we mentioned above, due to existence of  multiple closed geodesics
the obtained form of the heat kernel for compact symmetric spaces 
requires an additional regularization, which consists simply in an
analytical continuation of the result from the complexified noncompact
case. In any case, it gives a generating function for all heat
invariants and reproduces correctly the whole asymptotic expansion of
the heat kernel diagonal. However, since there are no closed geodesics
on non-compact symmetric spaces, it seems that the analytical
continuation of the obtained result for the heat kernel diagonal should
give the {\it exact result} for the non-compact case, and, even more
generally, for the general case too.  We have seen on the example of the
two-sphere that our method gives not just the asymptotic expansion of
the heat kernel diagonal but, after an appropriate regularization, in
fact, an exact result for the heat kernel diagonal.

%===================================================================


\begin{thebibliography}{99}

\bibitem{anderson90} 
A. Anderson and R. Camporesi, {\it Intertwining operators for solving 
differential equations with applications to symmetric spaces}, 
Commun. Math. Phys., {\bf 130} (1990), 61--82.

\bibitem{avramidi87}
I. G. Avramidi, {\it Covariant methods for the calculation of the
effective action in quantum field theory and investigation of
higher-derivative quantum gravity}, PhD Thesis, Moscow State University
(1987), arXiv:hep-th/9510140.

\bibitem{avramidi89} 
I. G. Avramidi, {\it Background field calculations in quantum 
field theory (vacuum polarization)}, 
Teor. Mat. Fiz., {\bf 79} (1989), 219--231. 

\bibitem{avramidi90} 
I. G. Avramidi, {\it The covariant technique for calculation of 
the heat kernel asymptotic expansion}, 
Phys. Lett. B, {\bf 238} (1990), 92--97.

\bibitem{avramidi91} 
I. G. Avramidi,  {\it A covariant technique for the calculation of the
one-loop effective  action}, Nucl. Phys. B {\bf 355} (1991) 712--754;
Erratum:  Nucl. Phys. B {\bf 509} (1998) 557-558. 

\bibitem{avramidi93} 
I. G. Avramidi, {\it A new algebraic approach for calculating 
the heat kernel in gauge theories}, 
Phys. Lett. B, {\bf 305} (1993), 27--34.

\bibitem{avramidi94a}
I. G. Avramidi,  {\it Covariant methods for calculating the low-energy
effective action in quantum field theory and quantum gravity}, 
University of Greifswald (March, 1994),  arXiv:gr-qc/9403036, 48 pp. 

\bibitem{avramidi94} 
I. G. Avramidi,  {\it The heat kernel on symmetric spaces via
integrating over the group  of isometries}, Phys. Lett. B {\bf 336}
(1994) 171--177.

\bibitem{avramidi95} 
I. G. Avramidi,  {\it Covariant algebraic method for calculation of the
low-energy heat kernel},  J. Math. Phys. {\bf 36} (1995) 5055--5070;
Erratum: J. Math. Phys. {\bf  39} (1998) 1720.

\bibitem{avramidi96} 
I. G. Avramidi,  {\it A new algebraic approach for calculating the heat
kernel in quantum  gravity}, J. Math. Phys. {\bf 37} (1996) 374--394.

\bibitem{avramidi98}
I. G. Avramidi, {\it Covariant approximation schemes  for calculation of
the heat kernel in quantum field theory}, in: {\it Quantum Gravity},
Eds. V. A. Berezin, V. A. Rubakov and D. V. Semikoz, World Scientific,
Singapore, 1998,  pp.~61--78.

\bibitem{avramidi99}  
I. G. Avramidi, {\it Covariant techniques for computation of the heat  
kernel}, Rev. Math. Phys., {\bf 11} (1999), 947--980. 

\bibitem{avramidi00} 
I. G. Avramidi, {\it Heat Kernel and Quantum Gravity}, Lecture Notes in
Physics, Series Monographs, LNP:m64, Springer-Verlag, Berlin, 2000.

\bibitem{avramidi02}
I. G. Avramidi, {\it Heat kernel approach in quantum field theory},
Nucl. Phys. Proc. Suppl., {\bf 104} (2002), {3--32}.

\bibitem{avramid06}
I. G. Avramidi, {\it Heat kernel asymptotics on symmetric spaces},
Int. J. Geom. Topol., (2007) (to be published); arXiv:math.DG/0605762 

\bibitem{barut77}
A. O. Barut and R. Raszka, {\it Theory of Group Representations and
Applications}, PWN, Warszawa, 1977.

\bibitem{berline92} 
N. Berline, E. Getzler and M. Vergne, {\it Heat Kernels and Dirac
Operators}, Springer-Verlag, Berlin, 1992.

\bibitem{camporesi90} 
R. Camporesi, {\it Harmonic analysis and propagators on homogeneous
spaces}, Phys. Rep. {\bf 196}, (1990), 1--134.

\bibitem{dowker70} 
J. S. Dowker, {\it When is the ``sum over classical paths'' exact?}, 
J. Phys. A, {\bf 3} (1970), 451--461.

\bibitem{dowker71} 
J. S. Dowker, {\it Quantum mechanics on group space and
Huygen's principle}, Ann. Phys. (USA), {\bf 62} (1971), 361--382.

\bibitem{erdelyi53} 
A. Erd\'elyi, W. Magnus, F. Oberhettinger and F. G. Tricomi,  {\it
Higher Transcendental Functions}, (McGraw-Hill, New York, 1953), vol. I.

\bibitem{fegan83}
H. D. Fegan, {\it The fundamental solution of the heat equation on a compact 
Lie group}, J. Diff. Geom., {\bf 18} (1983), 659--668.

\bibitem{gilkey75} 
P. B. Gilkey, {\it The spectral geometry of Riemannian manifold}, J.
Diff. Geom., {\bf 10} (1975), 601--618.

\bibitem{gilkey95}
P. B. Gilkey, {\it Invariance Theory, the Heat Equation and the
Atiyah-Singer Index Theorem}, CRC Press, Boca Raton, 1995.

\bibitem{helgason84}
S. Helgason, {\it Groups and Geometric Analysis: Integral Geometry,
Invariant Differential Operators, and Spherical Functions}, Mathematical
Surveys and Monographs, vol. 83, AMS, Providence, 2002, p. 270 

\bibitem{hurt83}
N. E. Hurt, {\it Geometric Quantization in Action: Applications of
Harmonic Analysis in Quantum Statistical Mechanics and Quantum Field
Theory}, D. Reidel Publishing, Dordrecht, Holland, 1983.

\bibitem{kirsten01}
K. Kirsten, {\it Spectral Functions in Mathematics and Physics}, 
CRC Press, Boca Raton, 2001.

\bibitem{ruse61}
H. Ruse, A. G. Walker and T. J. Willmore, {\it Harmonic Spaces},
Edizioni Cremonese, Roma, 1961.

\bibitem{takeuchi91} 
M. Takeuchi, {\it Lie Groups II}, in: {\it Translations of  Mathematical
Monographs}, vol. 85, AMS, Providence, 1991, p.167.

\bibitem{vandeven98}
A. E. M. Van de Ven,
{\it Index free heat kernel coefficients}, Class. Quant. Grav. {\bf 15}
(1998), 2311--2344. 

\bibitem{vassilevich03}
D. V. Vassilevich, {\it Heat kernel expansion: user's manual}, Phys.
Rep., {\bf 388} (2003), 279--360.

\bibitem{wolf72}
J. A. Wolf, {\it Spaces of Constant Curvature}, University of
California, Berkeley, 1972.

\bibitem{yajima04}
S. Yajima, Y. Higasida, K. Kawano, S.-I. Kubota, Y. Kamo and S. Tokuo,
{\it Higher coefficients in asymptotic expansion of the heat kernel},
Phys. Rep. Kumamoto Univ., {\bf 12} (2004), No 1, 39--62.

\end{thebibliography}
\end{document}